\documentclass{amsart}
\usepackage{amsmath}
\usepackage{amssymb}
\usepackage[all]{xy}
\usepackage{graphicx}
\usepackage{mathrsfs}
\usepackage{mathabx}
\usepackage{mathrsfs}

\numberwithin{equation}{section}

\usepackage{url}
\usepackage{longtable,tabu}
\usepackage{float}

\usepackage[latin1]{inputenc}
\usepackage{xspace,amssymb,amsfonts,euscript}
\usepackage{amsthm,amsmath}
\usepackage{palatino}
\usepackage{euscript}
\input xy \xyoption {all}

\usepackage{tikz,environ}
\usetikzlibrary{arrows}
\usetikzlibrary{patterns,snakes}

\RequirePackage{color}
\definecolor{myred}{rgb}{0.75,0,0}
\definecolor{mygreen}{rgb}{0,0.5,0}
\definecolor{myblue}{rgb}{0,0,0.65}

\RequirePackage{ifpdf}
\ifpdf
  \IfFileExists{pdfsync.sty}{\RequirePackage{pdfsync}}{}
  \RequirePackage[pdftex,
   colorlinks = true,
   urlcolor = myblue, 
   citecolor = mygreen, 
   linkcolor = myred, 
   pagebackref,
   bookmarksopen=true]{hyperref}
\else
  \RequirePackage[hypertex]{hyperref}
\fi

\RequirePackage{ae, aecompl, aeguill} 

    \def\CM{{\mathbb{C}}}
    \def\DM{{\mathbb{D}}}
    
    \def\FM{{\mathbb{F}}}

    \def\PM{{\mathbb{P}}}
    \def\QM{{\mathbb{Q}}}
    \def\RM{{\mathbb{R}}}

    \def\ZM{{\mathbb{Z}}}

    \def\AC{{\mathcal{A}}}

    \def\EC{{\mathcal{E}}}
    \def\FC{{\mathcal{F}}}
    
    \def\HC{{\mathcal{H}}}

    \def\laC{{\mathcal{L}}}

    \def\OC{{\mathcal{O}}}
    
    \def\QC{{\mathcal{Q}}}
    
    \def\SC{{\mathcal{S}}}

\def\ES{{\EuScript E}}
\def\FS{{\EuScript F}}
\def\GS{{\EuScript G}}

\newcommand{\nc}{\newcommand} \newcommand{\renc}{\renewcommand}

\def\a{\alpha}
\def\b{\beta}
\def\g{\gamma}

\def\e{\varepsilon}
\renc{\l}{\lambda}

\newcommand{\rdots}{\mathinner{ \mkern1mu\raise1pt\hbox{.}
    \mkern2mu\raise4pt\hbox{.}
    \mkern2mu\raise7pt\vbox{\kern7pt\hbox{.}}\mkern1mu}}

\def\sgn{{\mathrm{sgn}}}

\def\res{{\mathrm{res}}}

\DeclareMathOperator{\Irr}{Irr}

\def\pr{{\mathrm{pr}}}
\def\inc{{\mathrm{inc}}}

\def\un{\underline}

\def\p{{}^p}
\def\f{{}^f}

\def\to{\rightarrow}

\def\laongto{\laongrightarrow}

\def\onto{\twoheadrightarrow}
\nc{\triright}{\stackrel{[1]}{\to}}
\nc{\laongtriright}{\stackrel{[1]}{\laongto}}

\nc{\Hb}{H^\bullet}

\nc{\Br}{\mathcal{B}}
\nc{\HotRR}{{}_R\mathcal{K}_R}
\nc{\HotR}{\mathcal{K}_R}
\nc{\excise}[1]{}
\nc{\defect}{\text{df}}
\nc{\h}[1]{\underline{H}_{#1}}

\nc{\Ga}{\mathbb{G}_a} 
\nc{\Gm}{\mathbb{G}_m} 

\nc{\Perv}{{\mathbf{P}}}

\nc{\IH}{{\mathrm{IH}}}

\nc{\ic}{\mathbf{IC}}

\nc{\gl}{{\mathfrak{gl}}}
\renc{\sl}{{\mathfrak{sl}}}
\renc{\sp}{{\mathfrak{sp}}}

\renc{\Im}{\textrm{Im}}

\nc{\HBM}{H^{BM}}

\DeclareMathOperator{\For}{For} 
 \DeclareMathOperator{\Hom}{Hom}
 \DeclareMathOperator{\ch}{ch}

\DeclareMathOperator{\Rep}{\mathrm{Rep}}
\DeclareMathOperator{\Parity}{\mathrm{Par}}
\DeclareMathOperator{\Tilt}{\mathrm{Tilt}}

\DeclareMathOperator{\id}{id}

\newtheorem{thm}{Theorem}[section]

\newtheorem{cor}[thm]{Corollary}
  
\newtheorem{conj}[thm]{Conjecture}

\theoremstyle{definition}

\newtheorem{ex}[thm]{Example}
\newtheorem{warning}[thm]{Warning}

\theoremstyle{remark}
\newtheorem{remark}[thm]{Remark}

\DeclareMathOperator{\Ext}{Ext}

\newcommand{\ra}{\rightarrow}

\newcommand{\into}{\hookrightarrow}

\def\Gr{{\EuScript Gr^\vee}}
\def\Grx{{\EuScript Gr_x^\vee}}
\def\Grl{{\EuScript Gr_\lambda^\vee}}

\def\Fl{{\EuScript Fl^\vee}}
\def\Flx{{\EuScript Fl_x^\vee}}

\nc{\simto}{\stackrel{\sim}{\to}}

\DeclareMathOperator{\socle}{socle}

\DeclareMathOperator{\Lie}{Lie}

\DeclareMathOperator{\ind}{ind}

\nc{\ext}{\textrm{ext}}
\def\acts{\;\lefttorightarrow \;}
\nc{\fW}{{}^f W}
\nc{\pdot}{ \bullet_p}
\nc{\wdot}{ \bullet}
\nc{\la}{\langle}
\renc{\ra}{\rangle}

\nc{\Wf}{W}
\nc{\Wa}{\mathcal{W}}
\nc{\Sa}{\mathcal{S}}
\nc{\Wae}{\mathcal{W}^{\mathrm{ext}}}

\nc{\Fr}{\mathrm{Fr}} 
\nc{\Repp}{\Rep_0}
\nc{\Repe}{\mathrm{Rep}_0^{\mathrm{ext}}}

\nc{\AntiS}{\mathrm{AS}}
\nc{\CAS}{\mathcal{AS}}

\nc{\HCe}{\HC^{\mathrm{ext}}}

\nc{\Kar}{\mathrm{Kar}}
\nc{\bs}{\mathrm{BS}}

\nc{\Iw}{\mathrm{Iw}}
\nc{\ev}{\mathrm{ev}}

\newcommand{\bk}{\Bbbk}

\nc{\SL}{\mathrm{SL}}
\nc{\GL}{\mathrm{GL}}
\nc{\Sp}{\mathrm{Sp}}
\nc{\Eeight}{\mathrm{E}_8}
\nc{\Gtwo}{\mathrm{G}_2}

\nc{\He}{\mathrm{H}} 
\nc{\Hee}{\mathrm{H}^\mathrm{ext}} 
\nc{\Zvv}{\mathbb{Z}[v]}
\nc{\Zv}{\mathbb{Z}[v^{\pm 1}]}

\nc{\HCat}{\mathcal{H}} 
\nc{\HCate}{\mathcal{H}^\mathrm{ext}} 

\nc{\pcan}{{}^p\un{h}}

\title[Algebraic representations and constructible sheaves]{Algebraic representations and \\constructible
   sheaves}

\author{Geordie Williamson}
\address{Research Institute for Mathematical Sciences,
Kyoto University,
Kyoto 606-8502,
JAPAN.}

\date{\today}

\begin{document}

\begin{abstract} 
These are notes for my Takagi lecture at the University of Tokyo in
November, 2016. I survey what is known about simple modules for reductive algebraic
groups. The emphasis is on characteristic $p>0$ and Lusztig's character
formula. I explain ideas connecting representations and
constructible sheaves (Finkelberg-Mirkovi\'c conjecture) in the spirit
of the Kazhdan-Lusztig conjecture. I also discuss a conjecture
with S. Riche (a theorem for $\GL_n$) which should eventually make computations more feasible.
\end{abstract}

\maketitle

\section*{Introduction}

Let $G$ denote an algebraic group over an algebraically closed field
$\bk$. A representation of $G$ is a $\bk$-vector space $V$ and a
homomorphism $G \to \GL(V)$ of algebraic groups. In this article we discuss various approaches
to the representation theory of reductive algebraic groups (like
$\GL_n, \Sp_{2n}, \dots, \Eeight$) via constructible sheaves.

Studying the
representation theory of $G$ can be thought of as ``harmonic analysis
in algebraic geometry''. 
Over
fields of characteristic zero the theory is
well understood and extremely useful. It parallels the theory of compact
Lie groups. Much research over the last five decades has focused on the case of
characteristic $p >0$.
Here the theory is highly developed, however
several fundamental questions remain unsolved.

The deepest result in the field (at least on the level of characters)
is Lusztig's formula. It gives character formulas for certain simple
modules, from which the
characters of all simple modules can be deduced.\footnote{If our
  characteristic $p$ is not too small. Such subtleties will be ignored in the introduction.}
If we fix the root system of our group and let $p$ vary, then we know
that Lusztig's character formula holds if $p$ is very large. However only
in very few cases (e.g. $\SL_2, \SL_3, \SL_4, \Sp_{4}, \Gtwo$) do we know
precisely when it holds! We also don't understand well what happens
when it fails.

Lusztig's character formula was motivated by the
Kazhdan-Lusztig conjecture, which gives the characters of simple
highest weight representations of complex
semi-simple Lie algebras. The Kazhdan-Lusztig conjecture was first
proved by establishing a bridge to constructible sheaves on the flag
variety. Once one has traversed such a bridge, deep theorems concerning
constructible sheaves (e.g. the decomposition theorem, the Weil
conjectures, \dots) can be used to deduce the Kazhdan-Lusztig
conjecture.\footnote{In the words of Bernstein \cite{BernsteinNotes}: ``The amazing feature
  of the proof is that it does not try to solve the problem but just
  keeps translating it in languages of different areas of mathematics
  (further and further away from the original problem) until it runs
  into Deligne's method of weight filtrations which is capable to
  solve it.''}

By analogy with the Kazhdan-Lusztig conjecture one would like to build
a bridge between representations of $G$ and constructible sheaves. The
goal being to better understand Lusztig's character formula (amongst other things). Building such a bridge turns out to be much
harder in this setting. The most satisfactory such
statements are the geometric Satake equivalence and the
Finkelberg-Mirkovi\'c conjecture\footnote{The reader is warned
  that the Finkelberg-Mirkovi\'c conjecture is still a
  conjecture. However it is very useful as a
  guiding principle. Furthermore, recent work of Achar, Mautner,
  Riche and Rider
  seems to bring us close to a proof.}. Both results purport an equivalence
between the representation theory of $G$ and a category of perverse
sheaves on the affine Grassmannian $\Gr$ associated to the (complex) Langlands dual
group. Under both such equivalences the base field of the
representation theory corresponds to the coefficients of the perverse
sheaves. The space $\Gr$, however, is fixed.

The Finkelberg-Mirkovi\'c conjecture is easily seen to imply Lusztig's
character formula for large $p$. It also gives a character formula
for all $p$ in terms of the Euler characteristic of the stalks of
intersection cohomology complexes with $\bk$-coefficients in. In this way, deciding for which $p$ Lusztig's character formula
holds becomes a question about controlling torsion in certain local
integral intersection cohomology groups. Roughly speaking, it was by
producing many unexpected torsion classes that the author was recently
able to show that Lusztig's character formula cannot hold with the
hoped-for bounds.

The Finkelberg-Mirkovi\'c conjecture provides a very
conceptually satisfying ``constructible picture'' of representations. However it seems unlikely that it will help
with computations (at least with current tools). Indeed, the calculation of (the Euler
characteristics of) the stalks of intersection cohomology complexes
with coefficients in a field of characteristic $p > 0$ is notoriously
difficult.

In practice it is often easier to calculate the stalks of parity
sheaves. (These are analogues of intersection cohomology
complexes whose stalks satisfy a parity vanishing property. In this
setting they only
really become interesting with coefficients of positive characteristic.) Thus one
is led to try to find a character formula in which the stalks of
parity sheaves appear. Such a conjecture has recently been formulated
by Riche and the author, and proved for $G = \GL_n$. The result is a character formula
for tilting modules in terms of the $p$-canonical basis. This
conjecture should be related via Koszul duality to
the Finkelberg-Mirkovi\'c conjecture.

Due to limitations (both of time and the author's competence) we do not
discuss closely related categories of coherent sheaves. One can regard
algebraic representations
as $G$-equivariant coherent sheaves on a point. From this point of view most of the results
of this paper can be viewed as special cases of coherent / constructible
equivalences appearing in the geometric Langlands program. It
was in this context that characteristic zero analogues of the results we discuss
were often first proved \cite{ABG, AB, Baffine}. Another glaring
omission is that we do not discuss the infinitesimal group
schemes (Frobenius kernels etc.) which appear naturally in the
theory. Thus we do not discuss Lie algebra representations, nor the
Bezrukavnikov-Mirkovi\'c-Rumynin theory of localisation in positive
characteristic. This theory is
the natural extension to positive
characteristic of the original 
proof of the Kazhdan-Lusztig conjectures via $D$-modules.

\newpage

\subsection{Structure of this paper} This paper consists of two sections:
\begin{enumerate}
\item[\S \ref{sec:reps}.]\emph{Algebraic representations:} We review
  the fundamentals of the theory of representations of algebraic groups: classification
  of simple modules, induced modules, Weyl modules, tilting modules,
  Steinberg's tensor product theorem, the translation and linkage
  principles. Our goal is to give all results
  needed to understand the statement of Lusztig's character
  formula. We survey what is known and not known regarding
  Lusztig's formula. Finally, we explain an observation of
  Lusztig which predicts the values at 1 of certain
  affine Kazhdan-Lusztig polynomials.
\item[\S \ref{sec:sheaves}.]\emph{Constructible sheaves}: After reviewing the basics of
  perverse and parity sheaves we define the Hecke category. We then
  discuss the geometric Satake equivalence and Finkelberg-Mirkovi\'c
  conjecture. We explain why the Finkelberg-Mirkovi\'c conjecture
  implies character formulas for algebraic groups in terms of stalks
  of intersection cohomology complexes, and why the presence of
  torsion can be used to deduce that Lusztig's character formula does
  not hold for certain primes. Finally, we outline a conjectural link
  between  tilting modules and parity sheaves.
\end{enumerate}
We conclude the paper with a list of frequently used notation.

\subsection{Acknowledgements} I would like to thank P. Achar,
H. H. Andersen, R. Bezrukavnikov, C. Bonnaf\'e, S. Donkin, P. Fiebig, J. C. Jantzen, D. Juteau,
X. He, A. Henderson, S. Kumar, G. Lonergan, G. Lusztig, S. Makisumi, C. Mautner,
I. Mirkovi\'c,
S. Riche, L. Rider, R. Rouquier, W. Soergel and K. Vilonen for useful
discussions and observations on the subject of this paper. I am very
grateful to D. Juteau, M. Kaneda and S. Riche for feedback on
a first draft.

\subsection{Conventions}\label{sec:not}


If we write $G \acts X$ we mean that the group $G$ acts on 
$X$.

Given an abelian category $\AC$ we let $[\AC]$ denote the Grothendieck
group of $\AC$. If $\AC$ is an additive category, its split
Grothendieck group, denoted $[\AC]_\oplus$, is the quotient of the free
module on symbols $[M]$ for all objects $M \in \AC$ modulo the
relations $[M] = [M'] \oplus [M'']$ if $M \cong M' \oplus M''$. In
both settings the class of $M \in \AC$ is denoted $[M]$. If
$\AC$ is addition graded (i.e. equipped with an equivalence $M \mapsto
M[1]$) we view $[\AC]_\oplus$ as a $\ZM[v, v^{\pm
  1}]$-module via $v^{\pm 1}[M] := [M[\pm 1]]$.

\section{Algebraic representations} \label{sec:reps}

\subsection{Root data and the group}\label{sec:data}
We fix a reduced root datum $(X,
\Phi, X^\vee, \Phi^\vee)$ with $X$ the character lattice, $\Phi
\subset X$ the roots, $X^\vee$ the cocharacter lattice and $\Phi^\vee
\subset X^\vee$ the coroots. To our root datum we may associate a
split connected reductive ``Chevalley'' group
scheme $G_\ZM$ over $\ZM$. For any field $k$, extension of scalars
yields an algebraic group over $k$ which is split,
connected, reductive and has
the above root data. Throughout:
\begin{gather*}
  \text{$\bk$ denotes an algebraically closed field of characteristic $p
  \ge 0$}; \\
\text{$G_\bk$ denotes the reductive algebraic group over $\bk$ deduced from $G_\ZM$.}
\end{gather*}
We will make the following assumption:
\begin{equation*}
  \label{eq:4}
\begin{array}{c}  \text{Our root system (and thus our group $G_{\bk}$) is
    semi-simple} \\ \text{ and simply-connected: $\ZM \Phi^\vee = X^\vee$.}
\end{array}
\end{equation*}
(This assumption is not essential for most of the theory discussed
  below. However including it simplifies the exposition.) Let us fix a system of positive roots and coroots
\[
\Phi_+
\subset \Phi \quad \text{and} \quad \Phi_+^\vee \subset \Phi^\vee. \]
We let $T_\ZM \subset B_\ZM \subset G_\ZM$ denote 
``maximal torus'' and ``Borel'' subgroup schemes which arise in the
construction of $G_\ZM$; their extension of scalars yield a maximal torus
$T_\bk$ and Borel subgroup $B_\bk$ of $G_\bk$ for any $\bk$. We assume
that our Borel subgroup is chosen such that:
\begin{equation*}
\text{the positive roots $\Phi_+$ 
  are the weights that appear in $\Lie G / \Lie B$.}  
\end{equation*}
Thus the roots occurring in $\Lie B$ are the \emph{negative} roots
$-\Phi^+$. We denote by
\begin{gather*}
  X_+ := \{ \l \in X \; | \; \la \a^\vee, \l \rangle \ge 0 \text{ for
    all $\a^\vee \in \Phi^\vee_+$} \}, \\
  X^\vee_+ := \{ \g \in X^\vee \; | \; \la \g, \a \rangle \ge 0 \text{ for
    all $\a \in \Phi_+$} \}
\end{gather*}
the subsets of dominant weights and coweights.

\subsection{Representations and simple modules} \label{sec:simples} In
the next three sections we recall some fundamentals about
representations of reductive algebraic groups. The results are
standard and we do not give detailed references;
excellent sources include \cite{JantzenSurvey,  JaBook}.

Given a linear algebraic group $H$ defined over $\bk$ we denote by
$\Rep H$ its category of \emph{finite-dimensional} algebraic\footnote{It is traditional to call algebraic representations of $G_\bk$
  ``rational''. I am avoiding this
  terminology as it seems a reliable source of confusion for
  mathematicians from other fields.}
representations. What a representation is was defined in the introduction;
alternatively we could define $\Rep H$ to be the abelian category of
finite-dimensional $\bk[H]$-comodules, where $\bk[H]$ denotes the
regular functions on $H$. We denote by $\Irr H$ the set of isomorphism
classes of simple modules in $\Rep H$.

We will almost exclusively study representations of our semi-simple
group $G_\bk$. If the context is clear we will often abbreviate:
\[
\Rep := \Rep G_{\bk}.
\]

To any $B_{\bk}$-module $V$ we may associate the trivial vector
bundle $G_{\bk} \times V$ on $G_{\bk}$. The quotient for the $B_{\bk}$-action $b \cdot (g, v) :=
(gb^{-1}, bv)$ exists and yields a vector bundle $\laC_V$ on
$G_{\bk}/B_{\bk}$. Taking global sections of this vector bundle gives rise to the
induction functor:
\[
\ind_{B_{\bk}}^{G_{\bk}} : \Rep B_{\bk} \to \Rep G_{\bk} : V \mapsto \Gamma(G_{\bk}/B_{\bk}, \laC_V).
\]
This functor preserves finite-dimensional
modules because $G_\bk/B_\bk$ is complete.

In particular for any character $\l \in X$ of $T_\bk$ we can inflate
via $B_\bk \onto B_\bk/[B_\bk,B_\bk] = T_\bk$ to obtain a $B_\bk$-module
$\bk_\lambda$ and then induce (we set $\laC_\l := \laC_{\bk_\l}$)
\[
\ind_{B_{\bk}}^{G_{\bk}} \bk_\lambda= \Gamma(G_{\bk}/B_{\bk}, \laC_\lambda).
\]
It turns out that $\ind_{B_{\bk}}^{G_{\bk}} \bk_\lambda \ne 0$ if and
only if $\l \in X_+$. Thus, for $\l \in X_+$ we set
\[
\nabla_\l := \ind_{B_{\bk}}^{G_{\bk}} \bk_\lambda \in \Rep.
\]
We call $\nabla_\l$ an \emph{induced module}. If $p = 0$ then each
$\nabla_\l$ is simple. In general each $\nabla_\l$ has simple socle. We set
\[
L_\l := \socle(\nabla_\l).
\]
The following gives the classification of the simple $G_\bk$-modules:
\begin{thm} We have a bijection: \label{thm:class}
\begin{align*}
  X_+ & \simto \Irr G_\bk \\
  \l & \mapsto L_\l.
\end{align*}  
\end{thm}

We denote by $\sigma$ a Chevalley involution on $G_\bk$ and consider
the contravariant functor $\DM$ given by
\begin{gather*}
  V \mapsto (V^*)^{\sigma}
\end{gather*}
where $(-)^\sigma$ denotes twisting by the Chevalley involution. Then
$\DM$ is a duality on $\Rep$ (i.e. $\DM^2 \cong \id$). The twist by
$\sigma$ is to ensure
\begin{gather}
  \label{eq:DL}
\DM(L_\l) = L_\l.  
\end{gather}
We set
\[
\Delta_\l := \DM(\nabla_\lambda)
\]
and call it a \emph{Weyl module}. We could alternatively have defined
$L_\l$ as the simple head of $\Delta_\l$. For any $\l \in X_+$ we
have maps (unique up to scalar):
\begin{equation}
  \label{eq:delLnab}
\Delta_\l \onto L_\l \into \nabla_\l.  
\end{equation}

\begin{ex} \label{ex:SL21}
  If $G_\bk = \SL_2$ then we can identify $X = \ZM$, $\Phi = \{ \pm 2
  \}$, $X_+ = \ZM_{\ge
    0}$. We have $G_\bk/B_{\bk} = \PM^1$ and $\mathcal{L}_n = \OC(n)$ all $n \in
  X$. We have $\Gamma(\PM^1, \OC(n)) \ne 0$ if and only if $n \ge
  0$. If $V = \bk x \oplus \bk y$ denotes the natural module of
  $G_\bk$ then $\nabla_n = \Gamma(\PM^1, \OC(n)) = S^n(V)$ and
  $\Delta_n = \nabla_n^*$ for all $n \ge 0$. If
  $p = 0$ then all $\nabla_n$ are simple. If $p \ge 0$ then $\nabla_0,
  \dots, \nabla_{p-1}$ are simple but $\nabla_p$ is not: $L_p = \bk x^p \oplus
  \bk y^p \subset \nabla_p$ is a non-trivial submodule.
\end{ex}

\subsection{Characters} \label{sec:characters}
Any $M$ in $\Rep T_\bk$ is semi-simple and
$\Irr T_\bk = X$. Hence we have a canonical isomorphism
\begin{equation}
  \label{eq:repT}
[ \Rep T_\bk] = \ZM[X].  
\end{equation}
We identify both sides of \eqref{eq:repT} and write elements as
(finite) sums $\sum_{\l \in X} m_\l e^\l$. Given any $M \in \Rep
G_\bk$ its \emph{character}
\[
\ch M \in \ZM[X]
\]
is the class of the restriction of $M$ to
$T_\bk$ in $[\Rep T_\bk]$. Concretely,
\[
\ch M = \sum_{\l \in X} (\dim M(\l)) e^\l
\]
where $M(\l) \subset M$ denotes the $\l$ weight space of $T_\bk$.

Let $\Delta \subset \Phi^+$ denote the simple roots corresponding to
our choice of positive roots. Set
\[
\rho := \frac{1}{2} \sum_{\a \in \Phi_+} \a.
\]
Let $W$ denote the Weyl group with simple reflections $S = \{ s_\alpha
\; | \; \alpha \in \Delta \}$. We denote by $x \mapsto \e_x$ the sign
character of $W$. The \emph{dot action} of $W$ on $X$ is given by 
\[
x \wdot \l := x(\l + \rho) - \rho.
\]
For any $\l \in X^+$ consider the \emph{Weyl character}
\[
\chi_\l := \frac{\sum_{x \in W} \e_x e^{x \wdot \l} } { \sum_{x \in W} \e_x e^{ x
  \wdot 0}} \in \ZM[X]^W.
\]

If $p = 0$ then for any $\l \in X_+$ we have $\Delta_\l = L_\l =
\nabla_\l$ and
\[
\ch L_\l = \chi_\l.
\]
If $p > 0$ then this is no longer true in general, as we have already
seen for $SL_2$. 
However it is still true (a consequence of Kempf
vanishing: $H^i(G_\bk/B_\bk, \mathcal{L}_\l) = 0$ for $\l \in X_+$ and
$i > 0$) that
\begin{gather}
  \label{eq:Kempf}
\ch \Delta_\l = \ch \nabla_\l = \chi_\l.  
\end{gather}

The basic problem which
motivates this survey is:
\begin{equation}
  \label{eq:Q}
  \text{Determine $\ch L_\l$ for all $\l \in X_+$.}
\end{equation}
As explained above, the answer is known if $p = 0$. Thus $p > 0$ is
the case of interest for this survey. By considerations of highest weight 
\[
\{ [\Delta_\l] \; |\;  \l \in X_+ \}, \quad \{ [L_\l] \; |\;  \l \in X_+ \}
\quad \text{and} \quad \{ [\nabla_\l] \; | \l \in X_+ \} 
\]
are all bases for $[\Rep]$ (of course $[\Delta_\l] = [\nabla_\l]$ by
\eqref{eq:DL}). It turns out to be convenient to rephrase our basic
problem as follows:
\begin{equation}
  \label{eq:Q2}
  \text{Find expressions $[L_\l] = \sum_\mu m_{\mu,\lambda}[\Delta_\mu]$ for all $\l \in X_+$.}
\end{equation}
This is equivalent to writing $\ch L_\l$ in terms of Weyl
characters which, in turn, is equivalent to \eqref{eq:Q}.

\subsection{Steinberg's theorems} \label{sec:steinberg} Assume that $p
> 0$. Recall that our group $G_\bk$ arises by extension of scalars from a
group scheme $G_\ZM$ defined over the integers. In particular, it
arises via extension of scalars from a group over $\FM_p$ and hence has 
a natural $\FM_p$-rational structure. We denote by
\[
\Fr : G_\bk \to G_\bk
\]
the Frobenius map. (Concretely, because $G_\bk$ has an
$\FM_p$-rational structure, it can be defined as a closed subgroup of
some $GL_N$ by equations with coefficients in $\FM_p$; the Frobenius
map $\Fr$ is given by
the $p^{th}$-power on coordinates in any such embedding.) Precomposing by $\Fr$ defines the functor of \emph{Frobenius twist} on $\Rep$:
\[
M \mapsto M^{\Fr}.
\]
If $\ch M = \sum m_\l e^\l$ then $\ch (M^\Fr) = (\ch M)^\Fr := \sum m_\l e^{p
  \l}$. We denote the iterates of $\Fr$ by
\[
M \mapsto M^{\Fr^m}.\]

It is easy to see that if $M$ is simple, then so is
$M^\Fr$. (If $\bk$ is perfect then, as representations of abstract
groups, we are simply twisting by an automorphism.) However much more
is true. For any $\ell \ge 0$ consider the set of
\emph{$\ell$-restricted weights:}
\[
X_1^\ell = \{ \l \in X_+ \; | \; \langle \a^\vee, \l \rangle < \ell \text{
  for all $\a \in \Delta$} \}.
\]

\begin{thm}[Steinberg]
If $\l \in X_1^p$ and $\gamma \in X_+$ then $L_\l \otimes L_\gamma^\Fr$
is simple.
\end{thm}

By our assumption that our root system is simply connected there exist
fundamental weights $\{ \varpi_\a \; | \; \a \in \Delta \} \subset X_+$ (i.e. $\langle
\a^\vee, \varpi_\beta \rangle = \delta_{\a,\b}$ for all $\a, \b \in
\Delta$). We can rewrite $X_1^\ell$ in these coordinates as
\[
X_1^\ell = \{ \sum_{\a \in \Delta} a_\a \varpi_\a \; | \; 0 \le a_\a
< \ell \text{
  for all $\a \in \Delta$} \}.
\]
Given any $\l \in X_+$ we can consider its \emph{$p$-adic
  expansion}
\[
\l = \sum_{i = 1}^m \l_i p^i \quad \text{with $\l_i \in X_1^p$.}
\]
It follows immediately from Steinberg's theorem and Theorem
\ref{thm:class} that:
\begin{gather}
  \label{eq:stt}
 L_\l := L_{\l_0} \otimes L_{\l_1}^{\Fr} \otimes \dots \otimes
L_{\l_m}^{\Fr^m}.
\end{gather}

\begin{ex} We continue Example \ref{ex:SL21} with $G_\bk = \SL_2$. We
  have
\[
\ch \Delta_n = \ch \nabla_n = \chi_n = \frac{ e^{n} - e^{-n-2}}{e^0 - e^{-2}} =
e^n + e^{n-2} + \dots + e^{-n} \quad \text{for all $n \ge 0$.}
\]
Moreover, $\ch L_n = \ch \Delta_n$ if $n < p$ (i.e. if $n$ is
$p$-restricted). For general $n$ we consider its $p$-adic expansion $n
= \sum_{i = 0}^m n_ip^i$. By Steinberg's theorem:
\[
\ch L_n = (e^{n_0} + e^{n_0-2} + \dots + e^{-n_0})
(e^{n_1} + e^{n_1-2} + \dots + e^{-n_1})^\Fr \dots 
(e^{n_m} + \dots + e^{-n_m})^{\Fr^m}
\]
Thus Steinberg's theorem solves our basic question \eqref{eq:Q} for
$\SL_2$. However $\SL_2$ is essentially the only case
where Steinberg's theorem gives the complete answer.
\end{ex}

We now briefly recall the Steinberg restriction theorem. Logically it is 
irrelevant for the rest of this survey, however it is such a beautiful
theorem that it would be criminal not to mention it. Let us
temporarily denote by $G$ the split form of our group over
$\FM_p$. Everything that we have done in the previous sections can be done over
$\FM_p$. Hence we obtain representations $\Delta_\l, L_\l,
\nabla_\l$ for all $\l \in X_+$ of the group scheme $G$. Taking
rational points we obtain representations of the finite group of Lie
type $G(\FM_q)$ for any $q = p^\ell$. We denote these representations
by the same symbols.

\begin{thm}[Steinberg restriction theorem] 
The set $\{ L_\l \; | \;  \l \in X_1^q \}$ is a set of representatives
for the isomorphism classes of simple $\bk G(\FM_q)$-modules.
\end{thm}

In particular a solution to the basic question \eqref{eq:Q} would 
yield considerable information about the irreducible representations
of $\bk G(\FM_q)$ for all $q$.

\begin{remark}
  Steinberg's restriction theorem gives a remarkably tight connection
  on the level of simple modules. One further beautiful connection is given by the theory
  of generic cohomology \cite{CPSvdK, Parshall}. However on the level of
  categories the finite and algebraic groups appear quite
  different. At present we know much more
  about the category of algebraic representations than of $\bk
  G(\FM_q)$ (e.g. compare the induction theorems of \cite{ABG,
    HKS, ARFM} with the solution of Brou\'e's conjecture for
  $SL_2(\FM_{q})$ \cite{ChuangSL2,Okuyama}).
\end{remark}

\subsection{Tilting modules} \label{sec:tilting}
We briefly recall the theory of tilting
modules. Excellent sources for this material include the paper of
Donkin \cite{DonkinTilting}, the surveys of
Andersen \cite{ATilting} and Mathieu \cite{MatTilt} as well as \cite[Chapter E]{JaBook}.

The starting point is the fundamental
vanishing theorem:
\begin{gather}
  \label{eq:deltanabla}
\Ext^i(\Delta_\l, \nabla_\mu) = \begin{cases} \bk & \text{if $i = 0$ and
    $\l = \mu$,} \\
0 & \text{otherwise}. \end{cases}  
\end{gather}
Let us define $\Rep_\Delta$ (resp. $\Rep_\nabla$) to be the full
subcategory of $\Rep$ consisting of modules which admit a filtration
whose successive quotients are isomorphic to $\Delta_\mu$
(resp. $\nabla_\mu$) for some $\mu \in X_+$. We will call such a filtration a \emph{Weyl} (resp. \emph{good}) filtration.
 We say that a module is \emph{tilting} if it belongs to both $\Rep_\Delta$ and
$\Rep_\nabla$, that is, if it possesses both a Weyl and a good filtration. We denote by $\Tilt G_\bk$ (or $\Tilt$ if the context
is clear) the full subcategory of tilting modules. (Note that $\Tilt$
is additive but almost never abelian.)

Given a tilting module $M$ we denote by $(M:\Delta_\l)$ (resp. $(M:\nabla_\l)$ the 
multiplicity of $\Delta_\l$ in a Weyl (resp. good) filtration of
$M$. This number is well defined because
\[
\ch M = \sum_{\l \in X_+} (M:\Delta_\l) \chi_\l = \sum_{\l \in X_+} (M:\nabla_\l) \chi_\l.
\]

\begin{thm}[\cite{ringel}, \cite{DonkinTilting}] \label{thm:tiltclass}
  For each $\l \in X_+$ there exists an indecomposable tilting module $T_\l$ with
  highest weight $\l$. Moreover $\dim T_\l(\l) = 1$ and we have a
  bijection:
  \begin{align*}
    X_+ & \simto \left \{ \begin{array}{c} \text{indecomposable} \\
        \text{tilting modules} \end{array} \right \} / \cong
    \\
 \l & \mapsto T_\l.
  \end{align*}
\end{thm}

Note that $\DM$ exchanges $\Rep_\Delta$ and $\Rep_\nabla$ and thus preserves $\Tilt$. By highest weight considerations we deduce that
indecomposable tilting modules are self-dual: 
\begin{equation}
  \label{eq:tilt-self-dual}
\DM T_\l \cong T_\l.  
\end{equation}

If we write
\[
[T_\l] = \sum m_{\mu,\l} [\Delta_\mu]
\]
then $m_{\l,\l} = 1$ and $m_{\mu, \l} = 0$ if $\mu \not \le \l$ (again
by highest weight considerations). In particular, the elements
$[T_\l]$ are upper-triangular in the basis $\{ [\Delta_\l] \}$ of $[\Rep]$
and thus also provide a basis.

Another fundamental theorem concerning tilting modules is:

\begin{thm} \label{thm:ttt}
  \begin{enumerate}
  \item If $M, M'$ are tilting modules, then so is $M \otimes M'$.
  \item If $M$ is a tilting module and $L \subset G_\bk$ is a Levi
    subgroup, then the restriction of $M$ to $L$ is tilting.
  \end{enumerate}
\end{thm}

\begin{remark} 
The proof of Theorem \ref{thm:tiltclass} is not difficult. On the
other hand, Theorem \ref{thm:ttt} seems to be difficult. The first proof was given by Wang \cite{Wang} in type $A$ and large
characteristics for other groups, then Donkin \cite{DonkinLNM} gave a different
proof which covered almost all cases (he had to exclude $p = 2$ for
$E_7, E_8$). The first uniform proof is due to Mathieu and uses
Frobenius splitting \cite{MatFilt} (see also \cite{MatTilt}). For
other approaches to the theorem see
\cite{LittGood,PoloTilt,ParaTilt,KanedaTilt}. We will discuss an
approach to Theorem \ref{thm:ttt} via the affine Grassmannian in
\S\ref{sec:Satake}.
\end{remark}

\begin{remark}
Tilting modules provide powerful tools in the study of
$\Rep$.
\begin{enumerate}
\item Let $[\Tilt]_\oplus$ denote its split Grothendieck group. Theorem \ref{thm:ttt} implies that $\Tilt$ is a monoidal
  category, and thus $[\Tilt]_\oplus$ is a ring. The
  inclusion $\Tilt \into \Rep$ induces an
  isomorphism of rings\footnote{In
  the words of Donkin \cite{DonkinTilting}: ``Perhaps the main point of tilting modules is that they
provide a section of the character map.''}
\[
[\Tilt]_\oplus \simto [\Rep].
\]
Moreover, the classes $\{ [T_\l ] \; | \; \l \in X_+ \}$ give a basis with
strong positivity properties: it has positive coefficients when written in the basis
$[\Delta_\l]$; and, it has positive structure constants.
\item From \eqref{eq:deltanabla} it follows immediately that if $M, M'
  \in \Tilt$ then
\[
\Ext^i(M,M') = 0 \quad \text{for $i > 0$.}
\]
Moreover, $\Tilt$ is easily seen to generate the derived category
$D^b(\Rep)$. Tilting theory guarantees that the inclusion $\Tilt
\subset \Rep$ induces an equivalence of triangulated categories
\[
K^b(\Tilt) \simto D^b(\Rep).
\]
Thus tilting modules and morphisms between them provide a
``homological skeleton'' of $\Rep$. In this sense, tilting modules are
somewhat analogous to projective or injective
objects. (Note that neither injective nor projective objects exist in
$\Rep$: injective (resp. projective) objects only exist after passage
to the ind- (resp. pro-) completion of $\Rep$.)
\end{enumerate}
\end{remark}

\subsection{Tilting characters} \label{sec:tilting_char}
As well as the basic problem \eqref{eq:Q} of determining the character
of the simple modules, another problem which motivates this survey is:
\begin{equation}
  \label{eq:TQ}
  \text{Determine $\ch T_\l$ for all $\l \in X_+$.}
\end{equation}
As earlier it is convenient to reformulate the problem as follows:
\begin{equation}
  \label{eq:TQ2}
  \text{Find expressions $[T_\l] = \sum_\mu n_{\mu,\lambda}[\Delta_\mu]$ for all $\l \in X_+$.} 
\end{equation}
This appears to be a difficult problem. At the end of this
survey we will outline an approach to this problem via the Hecke
category. Here is a brief overview of what is known:
\begin{enumerate}
\item As for simple modules, there is a kind of tensor product theorem for
  tilting modules: if $M$ is indecomposable tilting and $\l$ belongs to the set
  $(p-1)\rho + X_1^p$ then
\[
T(\l) \otimes M^\Fr
\]
is indecomposable tilting. (More precisely, this is a theorem if $p
\ge 2h - 2$ by \cite[Proposition
2.1]{DonkinTilting} and Example 1 following it, and would follow for
all $p$ from \cite[Conjecture 2.2]{DonkinTilting}). This allows one
(under mild restrictions on $p$) to determine the
characters of all tilting modules from the knowledge of the characters
of $T(\l)$ belonging to the two sets
\begin{gather*}
  \{ \l \; | \; \l \in (p-1)\rho + X_1^p \} \quad \text{and} \\
  \{ \l = \sum_{\a \in \Delta} a_\a\varpi_\a \; | \; 0 \le a_\a < p \text{ for some $\a
    \in \Delta$}\}.
\end{gather*}
Note, however that the second set is infinite in all types other
than (products of) $\SL_2$. For $\SL_2$ this formula can be used to
determine all tilting characters, see \cite[\S 2, Example 2]{DonkinTilting}. See
\cite{LusztigWilliamson} for a description of the characters that
may be obtained in this way for a general root system.
\item By an observation of Andersen \cite{ATiltAlg} (see also \cite[\S
  1.8]{RW}), knowledge of a finite set
  of tilting characters determines all simple characters if $p \ge
  2h-2$.
\item By a result of Erdmann \cite{Erdmann} the determination of the
  characters of indecomposable tilting modules for $\GL_n$ (or
  equivalently $\SL_n$) in characteristic $p$ is equivalent to
  determining all decomposition numbers for representations of all
  symmetric groups which are indexed by partitions with less than or
  equal to $n$ parts. This is an unsolved problem for $n \ge
  3$, reflecting the fact that the characters of indecomposable tilting
  modules are unknown for $\SL_n$ for $n \ge 3$.
\item By results of Donkin on the Ringel self-duality of the Schur
  algebra one can rephrase the question of determining the characters
  of a fixed tilting module for $\GL_n$ in terms of certain simple
  characters for $\GL_N$ for some (usually much larger) $N$. Thus
  knowledge of all tilting characters for some $\GL_n$ and fixed
  characteristic $p$ would yield some information about simple
  characters for $\GL_N$ for large $N$ (and thus for $p$ ``small'' relative to
  $N$).
\item Andersen \cite{AFilt1,AFilt2} and Andersen-Kulkarni
  \cite{AndersenKulkarni} have proved a sum formula for tilting
  modules (this formula was inspired by Jantzen's sum formula for Weyl
  modules). Like Jantzen's formula it does not give complete
  information, but is very useful in small rank. Jensen \cite{JensenTilting} (see also
  Parker \cite{Parker}) has used this formula to determine some new tilting
  characters for $\SL_3$.
\item In the analogous setting of quantum groups at a root of unity
  the determination of the characters of the indecomposable tilting
  modules was solved by Soergel \cite{SoeKL,SoergelKippKM}.
\end{enumerate}

\subsection{The (extended) affine Weyl group} \label{sec:Wa} Here we briefly
discuss the (extended) affine Weyl group.   A very clear treatment of
this material can be found in \cite{IM}.

Let $X_\RM := X \otimes_{\ZM} \RM$. The affine Weyl
group $\Wa$ is the subgroup of affine transformations of $X_\RM$
generated by $\Wf$ (acting linearly) and $\ZM \Phi$ (acting by
translation). In formulas:
\[
  \Wa := \Wf \ltimes \ZM \Phi \acts X_\RM.
\]
Given $\l \in \ZM \Phi$ we denote by $t_\l \in \Wa$ the corresponding
translation.

The group $\Wa$ is also generated by the affine reflections
\[
s_{\a, m}(\l) := \l - \langle \a^\vee, \l \ra \a + m\a
\]
in the hyperplanes
\[
H_{\a, m} := \{ \l \in X_\RM \; | \; \langle \a^\vee, \l \ra = m \}
\]
for all $\a \in \Phi^+$ and $m \in \ZM$. The set
\[
C_- := \{ \l \in X_\RM \; | \; -1 \le \la \a^\vee, \l \ra \le 0 \;
\text{for all $\a \in \Phi^+$} \; \} \subset X_\RM
\]
is a fundamental domain for the action of $\Wa$ on $X_\RM$
\cite[Ch. V, \S 3]{Bourbaki}.


Consider the set $\Sa$ of reflections in those hyperplanes $H_{\a, m}$
which intersect $C_-$ in codimension one (the \emph{walls} of $C_-$). Then
$\Sa$ generates $\Wa$. Moreover, equipped with these generators $\Wa$
is a Coxeter group \cite[Ch. V, \S 3]{Bourbaki}. Throughout, whenever we view $\Wa$ as a Coxeter
group, it will always be with respect to the generators $\Sa$. We
denote by $\ell$ the length function on $\Wa$ with respect to the
generating set $\Sa$ and by $\le$ the Bruhat order on $\Wa$.

\begin{warning}
  Most authors (for example \cite{Bourbaki, IM}) define the affine Weyl
  group to be the semi-direct product of $W$ with the \emph{co}root
  lattice. Thus, our $\Wa$ is what is usually referred to as the
  affine Weyl group of the dual root system $\Phi^\vee$. In
  particular, to determine the Coxeter type of $\Wa$, one should
  consider the extended Dynkin diagram of the \emph{dual} root system
  $\Phi^\vee$. The convention we adopt here is better adapted to the
  combinatorics of representations of algebraic groups. It can also be
  seen as a shadow of Langlands duality, as should become clearer in
  \S\ref{sec:sheaves}.
\end{warning}

The \emph{extended affine Weyl group} is the subgroup of affine
transformations of $X_\RM$ generated by $W$ and
the weight lattice $X$:
\begin{gather*}
 \Wae:= W \ltimes X \acts X_\RM.
\end{gather*}
As above, given $\l \in X$ we denote by $t_\l \in \Wae$ the
translation by $\l$. The extended affine Weyl group has a length
function $  \ell : \Wae \to \ZM_{\ge 0}$ given by
\begin{align*}
  \ell(x) = \left | \left \{ \begin{array}{c} \text{hyperplanes
        $H_{\a,n}$ separating the } \\ \text{interior of $C_-$ from
        that of $x(C_-)$} \end{array} \right \} \right |.
\end{align*}
Because $\ZM \Phi \subset X$, the affine Weyl group $\Wa$ is a
subgroup of $\Wae$. The length function $\ell$ restricts to the standard (Coxeter)
length function on $\Wa$.

If we consider the subset of \emph{length
  zero elements}
\[
\Omega := \{ \omega
\in \Wae \; | \; \ell(\omega) = 0 \} = \{ \omega \in \Wae \; | \; \omega (C_-) = C_- \}
\]
then $\Omega \cong X/ \ZM \Phi$, $\Omega$ acts via conjugation as automorphisms of the Coxeter
system $(\Wa, \Sa)$, and $\Wae$ is the semi-direct product
\[
\Wae = \Omega \ltimes \Wa.
\]
We extend the Bruhat order to $\Wae$ by declaring that $\omega_1x_1 \le
\omega_2 x_2$ for $\omega_1, \omega_2 \in \Omega$, $x_1, x_2 \in \Wa$
if $\omega_1 = \omega_2$ and $x_1 \le x_2$.

Of course, both $\Wa$ and $\Wae$ preserve the character lattice $X \subset
X_\RM$.



\subsection{The (extended) affine Hecke algebra and Kazhdan-Lusztig
  basis} \label{sec:hecke}
The \emph{extended affine
Hecke algebra} $\Hee$  is the $\Zv$-algebra generated by symbols $\{
h_w \; | \; w \in \Wae \}$ subject to the relations
\begin{gather*}
  h_wh_{w'} = h_{ww'} \quad \text{if $\ell(ww') = \ell(w) + \ell(w')$,
  and}\\
  h_s^2 = (v^{-1} - v)h_s + h_{\id} \quad \text{for $s \in \Sa$}.
\end{gather*}
It is an associative unital (with unit $1 := h_{\id}$) algebra. The
set $\{ h_x \; | \; x \in \Wae \}$ is a $\Zv$-basis for $\Hee$ called
the standard basis. The basis elements $\{ h_x \; | \; x \in \Wa \}$
(resp. $\{ h_x \; | \; x \in W \}$) span a subalgebra
$\He \subset \Hee$, the \emph{affine Hecke algebra}
(resp. \emph{finite Hecke algebra}). The affine Hecke
algebra is isomorphic to the Hecke
algebra of the Coxeter system $(\Wa, \Sa)$.

Each standard basis element $h_x$ is invertible. The Kazhdan-Lusztig
involution is the algebra involution $h \mapsto \overline{h}$ on
$\Hee$ determined by $h_x \mapsto h_{x^{-1}}^{-1}$ and $v \mapsto
v^{-1}$. The following is a 
classical theorem of Kazhdan-Lusztig:

\begin{thm}[Kazhdan-Lusztig \cite{KLp, Lq}]
  For all $x \in \Wae$ there exists a unique element $\un{h}_x$ such
  that:
  \begin{enumerate}
  \item (``self-duality'') $\overline{\un{h}_x} = \un{h}_x$;
  \item (``Bruhat upper-triangularity'') we have
\[
\un{h}_x = \sum_{y \le x} h_{y,x} h_y
\]
for polynomials $h_{y,x} \in \Zvv$ with $h_{x,x} = 1$ and
$h_{y,x} \in v\ZM[v]$ for all $y < x$.
  \end{enumerate}
\end{thm}

By property (2) the set $\{ \un{h}_x \}$ is a basis for $\Hee$, the \emph{Kazhdan-Lusztig basis}. The polynomials $h_{y,x}$ are the
\emph{Kazhdan-Lusztig polynomials}. We extend their definition to all pairs
$y, x$ by setting $h_{y,x} := 0$ if $y \not \le x$.

\begin{remark}
  For an excellent introduction to Kazhdan-Lusztig polynomials the
  reader is referred to \cite{SoeKL}.
\end{remark}

\begin{remark}
  For $\omega \in \Omega$ and $x \in \Wa$ we have
  $h_{\omega^{-1}} \un{h}_{\omega
    x} = \un{h}_{x}$ and hence $h_{\omega y, \omega x} = h_{
    y,x}$ for all $y \in \Wae$. Thus all Kazhdan-Lusztig polynomials may be
  calculated in $\He$ and have non-negative coefficients (see Theorem
  \ref{thm:klpoly}).
\end{remark}

\subsection{The linkage principle} \label{sec:linkage}

From now on we assume that $p > 0$. In classical highest weight representation theory it is
often necessary to shift the origin to $-\rho$ and consider the \emph{dot
action} of $W$ on $X_\RM^*$:
\[
x \wdot \mu := x(\mu + \rho) - \rho.
\]
In the representation theory
of $G_{\bk}$ in characteristic $p$ it is necessary to dilate the action of
the (extended) affine Weyl group by $p$ and shift the origin to
$-\rho$. In this way we are led to the \emph{$p$-dilated dot action}:
\[
\Wa \pdot \acts X_\RM^*  \quad ( \text{resp. }
\Wae \pdot \acts X_\RM^* \;)
\]
defined via
\begin{align*}
  x \pdot \mu &:= x \wdot \mu \quad \text{for $x \in \Wf$ and}\\
  t_\l \pdot \mu &:= \mu + p\l \quad \text{for $\l \in \ZM \Phi$ (resp. $\l
    \in X$).}
\end{align*}
We say that $\l \in X_\RM^*$ is \emph{$p$-regular} if its stabiliser
under the $p$-dilated dot action is trivial. After $p$-dilation and $-\rho$ shift, a fundamental domain for the
$\pdot$-action is
\[
C_-^p := \{ \mu \in X_\RM^* \; | \; -p \le \la \a^\vee , \mu + \rho \ra
\le 0 \; \text{for all $\a \in \Phi^+$} \; \}.
\]

For $\mu \in C_-^p$, let $\Rep_\mu$ denote the full 
subcategory of all algebraic representations of $G_{\bk}$ whose composition
factors are simple modules with highest weight belonging to the $\Wa$ orbit of
$\mu$ under the $p$-dilated dot action:
\[
\Rep_\mu := \la L_\l \; | \; \l \in W \pdot \mu \cap X_+
\ra \subset \Rep.
\]

The linkage principle asserts that we have a direct sum decomposition
of abelian categories:
\begin{equation}
  \label{eq:block_almost}
\Rep = \bigoplus_{\mu \in C_-^p} \Rep_\mu.  
\end{equation}
In other words, any indecomposable module belongs to some $\Rep_\mu$
and if $M \in \Rep_\mu, M' \in \Rep_{\mu'}$ then $\Hom(M,M') = 0$
unless $\mu = \mu'$. Abusing language, we will refer to each
$\Rep_\mu$ as a \emph{block} of $\Rep$.

\begin{remark} This is an abuse of language because the decomposition
  \eqref{eq:block_almost} is not the finest possible, and hence does
  not give the block decomposition in the usual meaning of the term.  To understand the true block
  decomposition one needs to consider analogues of the $p$-dilated dot
  action
  for higher powers of $p$ \cite{DonkinBlocks}. However below we will
  assume that $p$ is 
  greater than the Coxeter number in which case $\Rep_\mu$ is
  indecomposable as an abelian category for ``most'' $\mu \in C_-^p$.\footnote{More precisely, $\Rep_\mu$ is
    indecomposable if and only if there exists $\a \in \Phi_+$ such
    that $\la \a^\vee, \mu + \rho \rangle$ is not divisible by $p$
    \cite[\S II.7.2]{JaBook}.} 
\end{remark}

Consider the set
\[
C_+^p := \{ \l \in X_+ \; | \; \la \alpha^\vee, \mu + \rho \rangle
\le p\;  \text{for all $\alpha \in \Phi_+$}\}.
\]
It is a consequence of the linkage principle that
\begin{equation}
  \label{eq:delta_simple}
\Delta_\l \text{ is simple, if }  \l \in C_+^p.
\end{equation}

\subsection{The (extended) principal block} \label{sec:repp}

Of particular importance is the \emph{principal block}
\[
\Repp := \la L_\l \; | \; \l \in \Wa \pdot 0 \cap X_+
\ra
\]
and the \emph{extended principal block}
\[
\Repe := \la L_\l \; | \; \l \in \Wae \pdot 0 \cap X_+
\ra.
\]
Note that the trivial module $L_0$ belongs to $\Rep_0$ (which explains
the name principal block). Because $0 = w_0 \pdot (-2\rho)$, in the above notation we have
\[
\Repp = \Rep_{-2\rho}.
\]
Similarly, if $\Omega \subset \Wae$
denotes the subset of length zero elements defined in the previous section one has
\[
\Repe = \bigoplus_{\omega \in \Omega} \Rep_{\omega \pdot (-2\rho)}.
\]

The importance of the extended principal block is the
following. Consider the Frobenius twist functor
\[
(-)^{\Fr}: \Rep \to \Rep.
\]
It sends a simple module of highest weight $\l$ to a simple module of
highest weight $p\l$. In particular its image lands in the extended
principal block, and we may view Frobenius twist as a functor:
\[
(-)^{\Fr}: \Rep \to \Repe.
\]

In fact we can say a little more. For any weight $\l \in X$ we
can write $\l = \l_0 + p\l_1$ with $\l_0 \in X_1^p$ and
\[
L_\l = L_{\l_0} \otimes L_{\l_1}^\Fr.
\]
Hence
\[
L_{\l} \otimes L_{\g}^\Fr = L_{\lambda_0} \otimes ( L_{\lambda_1} \otimes
L_{\gamma})^\Fr.
\]
Thus if $L_{\l} \in \Repp$ then so is $L_{\l}
\otimes L_{\gamma}^\Fr \in \Repe$. In other words, the bifunctor
\begin{align*}
  \Repe \times \Rep & \to \Repe \\
  (V, M) & \mapsto V \otimes M^\Fr
\end{align*}
makes $\Repe$ a (right) module category over
$\Rep$.

\begin{remark}
  More generally, any ``extended'' block
\[
\Rep_\mu^\ext := \la L_\l \; | \; \l \in \Wae \pdot \mu \cap X_+
\ra
\]
is a module category over $\Rep$ via the Frobenius twist.
\end{remark}


\subsection{Translation functors} \label{sec:translation}
The linkage principle implies that in order to understand $\Rep$ as an abelian category it is
enough to understand each block $\Rep_\mu$. Translation functors can be used
to relate these blocks and often reduce questions to the study of the
principal block $\Rep_0$.

For any $\l \in C_-^p$ let $\inc_\l$ (resp. $\pr_\l$) denote the
inclusion (resp. projection) functor to the block $\Rep_\l \subset \Rep$. Fix $\l, \mu \in C_-^p$. We define the \emph{translation functor}
\[
T_\l^\mu : \Rep_\l \to \Rep_\mu
\]
via
\[
T_\l^\mu := \pr_\mu( V \otimes ( \inc_\l( - ) ))
\]
where $V$ is any module whose extremal weights are $W( \mu -
\l)$.\footnote{Different choices of module yield isomorphic functors
  \cite[\S 7.6, Remark 1]{JaBook}. Thus with the definition above $T_\l^\mu$ is
  only defined up to isomorphism.} For example, if $\nu$ is the unique element of $W ( \mu - \l)
\cap X_+$ we could take $V := L_\nu$ or $V := \Delta_\nu$. One may
think of a translation functor as a certain ``matrix coefficient of the functor
$V \otimes (-)$''. Because $V \otimes (-)$ and $V^* \otimes (-)$ are biadjoint
one easily deduces that $T_\l^\mu$ and $T_\mu^\l$ are biadjoint.

To describe the effect of translation functors on blocks we need a
little more notation. For any $x \in C_-$ its stabiliser in $\Wa$ is generated by those $s
\in \Sa$ which fix it. A \emph{facet} of $C_-$ is a non-empty subset
consisting of all points with a fixed stabiliser in $\Wa$. Each facet
is locally closed (i.e. open in its closure) and $C_-$ is the union of
its facets. The same statements and definitions apply verbatim for
$C_-^p$ if instead we consider the $p$-dilated dot action of $\Wa$.

The two most useful properties of translation functors, often called
\emph{Jantzen's translation principles}, are the
following:
\begin{enumerate}
\item If $\l, \mu \in C_-^p$ belong to the same facet then
$T_\l^\mu : \Rep_\l \to \Rep_\mu$
is an equivalence of abelian categories \cite[Proposition II.7.8]{JaBook}
preserving standard modules.
\item If $\l, \mu \in C_-^p$ and $\mu$ belongs to the closure of the
  facet containing $\l$ then $T_\l^\mu : \Rep_\l \to \Rep_\mu$ sends
  each simple (resp. standard) module to a simple (resp. standard)
  module or zero. We refer the reader to \cite[\S\S II.7.11-15]{JaBook}
  for the precise statements.
\end{enumerate}

 A consequence of these two properties is that, if $\mu$ belongs
  to the closure of the facet containing $\l$, and if one knows
  character formulas
  \begin{equation}
    \label{eq:char}
\ch L_{x \pdot \l} = \sum a_{y,x} \chi_{y \pdot \l}    
  \end{equation}
for all $x \pdot \l \in X_+$, then one may easily deduce similar
character formulas for $L_{x \pdot \mu}$ for all $x \pdot \mu \in
X_+$.

The interior of $C_-^p$ constitutes the unique open facet of
$C_-^p$. The following are equivalent:
\begin{enumerate}
\item there exists a $p$-regular weight $\l \in X$;
\item $C_-^p$ contains a point of $X$ in its interior;
\item $C_-^p$ contains $-2\rho$ in its interior;
\item $p \ge h$ where $h$ is the \emph{Coxeter
    number}\footnote{Warning: $h$ will usually disagree with the classical
    definition of the Coxeter number if our root system is decomposable.}:
\[
h = \max_{\a \in \Phi_+} ( \la \a^\vee, \rho \ra + 1).
\]
\end{enumerate}
Thus, if $p \ge h$ and we know expressions \eqref{eq:char} for all
simple modules in $\Repp$, then we may deduce character formulas for
all simple modules in $\Rep$.

\subsection{Lusztig's character formula} \label{sec:lcf}
We keep the notation from
previous sections.

\begin{conj}[Lusztig conjecture \cite{L80}, original version] \label{LC:original}
Fix a
  $p$-regular weight $\mu \in C_-^p$ and $x \in \Wa$ such that $x
  \pdot \mu \in X_+$. Suppose that $p \ge h$ and that
  $\la \a^\vee, x \pdot \mu + \rho \ra \le p(p-h+2)$ for all $\a \in
  \Phi_+$ (``Jantzen's condition''). Then
\begin{equation}
  \label{LCF} \tag{LCF}
  [ L_{x \pdot \mu} ] = \sum_{y \le x \atop y \pdot \mu \in X^+}
  \e_{yx} h_{y,x}(1) [ \Delta_{y \pdot \mu} ].
\end{equation}
\end{conj}

\begin{remark} \label{rem:LCF} Some remarks concerning Lusztig's conjecture:
  \begin{enumerate}
  \item \label{rem:LCFextended} Lusztig's original formulation fixed the choice $\mu = -2\rho$. It is
    equivalent to the above formulation by Jantzen's translation
    principle. In this way one can also see that \eqref{LCF} implies a similar
    formula where $x$
    and $y$ are allowed to belong to the extended affine Weyl group
    $\Wae$ rather than $\Wa$ (but still satisfy the other conditions).
  \item \label{rem:singular} If $\mu$ is not $p$-regular then using translation functors
    one can deduce from \eqref{LCF} an identical expression for $[
    L_{x \pdot \mu} ]$ as long as one assumes that $x$ is of minimal
    length amongst all such $x' \in \Wa$ with $x \pdot \mu = x' \pdot
    \mu$.
  \item One
    of the remarkable aspects of \eqref{LCF} is that it predicts that
    part of the representation theory of $G_{\bk}$ is 
    ``independent of $p$'': if we use the $p$-dilated dot
    action to parametrise our highest weights, the coefficients expression
    simple modules in terms of standard modules are independent of
    $p$!
\item \label{rem:JR} Let us try to explain the meaning of Jantzen's condition. 
Consider the $p$-adic expansion of
our highest weight
\[
x \pdot \mu = \sum \l_i p^i \quad \text{with $\l_i \in X_1^p$.}
\]
Jantzen noticed (see \cite[\S 4.4]{JaCF} and \cite[\S 8.22]{JaCF}) that a necessary condition for independence of
$p$ is that
\begin{equation}
  \label{eq:JR}
\text{$\l_i$ is zero for $i \ge 2$ and $\l_1 \in C_p^+$. }
\end{equation}
It is easy to prove that if $x \pdot \mu$ satisfies
Jantzen's condition then it satisfies \eqref{eq:JR}. However 
in general there will be weights satisfying \eqref{eq:JR} which do not
satisfy Jantzen's condition. As far as we can tell, Jantzen's
condition provides an easily defined and large set on which
\eqref{eq:JR} holds, and has no significance beyond that.
\item Let us try to give a rough idea why \eqref{eq:JR} is necessary
  for independence of $p$. Write
\[
x \pdot \mu = \l_0 + p\l_1
\]
with $\l_0 \in X_1^p$. Then $\l_1$ is independent of $p$. Consider the module
\[
\widetilde{L_{x \pdot \mu}} := L_{\l_0} \otimes \Delta_{\l_1}^\Fr.
\]
It follows from Lusztig's conjecture for quantum groups (a theorem)
that
\begin{equation*}
  [ \widetilde{L_{x \pdot \mu}} ] = \sum_{y \le x \atop y \pdot \mu \in X^+}
  \e_{yx} h_{y,x}(1) [ \Delta_{y \pdot \mu} ].
\end{equation*}
(The important point is that this expression is independent of $p$.)
If $p$ is large enough then $\l_1$ belongs to $C_p^+$ and $\Delta_{\l_1}$ is
simple, and thus we have equality $\widetilde{L_{x \pdot \mu}} =
L_{x \pdot \mu}$ by Steinberg's tensor product theorem. In particular,
if we can write $[L_{x \pdot \mu}]$ in terms of standard modules in a
manner which is independent of $p$ then we must have $\widetilde{L_{x
    \pdot \mu}} = L_{x \pdot \mu}$ and hence $\Delta_{\l_1}$ must be
simple. This is ensured by \eqref{eq:JR}.
\item   For other discussions of Lusztig's conjecture see
  \cite{ASurvey}, \cite{DLC}, \cite{Scott} and \cite{Humphreys}.
\end{enumerate}
\end{remark}

\subsection{Lusztig's conjecture and the Steinberg tensor
  product theorem} \label{sec:LCSteinberg}
Suppose that $\mu$ is $p$-regular and $x \pdot \mu$ is
dominant. Consider the $p$-adic expansion of
our highest weight
\[
x \pdot \mu = \sum \l_i p^i \quad \text{with $\l_i \in X_1^p$.}
\]
If we can apply \eqref{LCF} to $L_{x \pdot \mu}$
then $\l_i = 0$ for $i \ge 2$ and $\l_1 \in C_p^+$ (see Remark \ref{rem:LCF}\eqref{rem:JR}). Hence,
$\Delta_{\l_1}$ is simple and by Steinberg's tensor product theorem:
\[
L_{x \pdot \mu} = L_{\l_0} \otimes \Delta_{\l_1}^\Fr.
\]
If $\l_1 \ne 0$ (and assuming \eqref{LCF}) there are two ways to
calculate the character of this module:
\begin{enumerate}
\item We can apply \eqref{LCF} to $L_{x \pdot \mu}$ directly;
\item We can apply \eqref{LCF} to calculate $L_{\l_0}$ and then
  multiply it with the Frobenius twist of the character of
  $\Delta_{\l_1}$. (Note that $\l_0$ might no longer be $p$-regular, in
  which case we apply Remark \ref{rem:LCF}\eqref{rem:singular}.)
\end{enumerate}
Thus the following is reassuring:

\begin{thm} [Kato \cite{Kato}] Lusztig's conjecture is consistent with
  the Steinberg tensor product theorem; that is, both of the above approaches give the same answer.
\end{thm}

For any fixed $p$, Lusztig's conjecture provides us with only finitely
many characters. However, by Steinberg's tensor product theorem we can
calculate the characters of all $L_\l$ for $\l \in X_+$ if we know the
characters of $L_\l$ for all $\l \in X_1^p$. All weights in $X_1^p$
satisfy the Jantzen condition if and only if
\[
\la \a^\vee, (p-1)\rho + \rho  \ra \le p(p-h+2) \quad \text{for all
  $\a \in \Phi_+$,}
\]
or in other words if
\[
h-1 \le p - h + 2 \; \Leftrightarrow \; p \ge 2h-3.
\]
Thus Lusztig's conjecture provides a complete conjectural answer if $p \ge 2h-3$.

Since Kato's result several authors came to regard the following
stronger version of Lusztig's conjecture as realistic:\footnote{Scott
  \cite{Scott} refers to it as ``Kato's extension of the Lusztig
  conjecture''. Jantzen \cite{JaBook, JaCF} says that it ``seems to be
  a realistic conjecture''.}

\begin{conj}[Lusztig conjecture, revised version] \label{LC:revised}
Suppose that $p \ge
  h$. Then \eqref{LCF} holds if $x \pdot \mu \in X_1^p$.
\end{conj}

\begin{remark} \label{rem:affine_finite}
As we will discuss below, the bound both in the original and revised
version of Lusztig's conjecture are much too optimistic. However 
Kato's theorem and the revised version are important for (at least)
the following reason: as $p$ varies, the number of weights for which
one needs to check the original formulation of Lusztig's conjecture
gets larger and larger. However it is not difficult to see that the
set
\[
\{ w \in \Wa \; | \; w \pdot (-2\rho) \in X_1^p \}
\]
is independent of $p$, as long as $p \ge h$. In this way Lusztig's
conjecture becomes a finite problem and one might hope to settle it
for ``all primes at once''. We will have more to say about this in the
second part.
\end{remark}

\begin{remark}
  As we have mentioned above, combining Steinberg's tensor product
  theorem and \eqref{LCF} yields a character formula for all highest
  weights $\l \in X^+$. However for general $\l$ this is rather
  indirect: one needs to apply \eqref{LCF} once for each $p$-adic
  digit of $\l$. Haboush \cite{Haboush}  and Humphreys \cite[\S
  3.12]{Humphreys}  advocate the consideration of different Weyl groups for each power
  of $p$. This intriguing idea appears not yet to have borne fruit. \end{remark}

\begin{remark}
Recently Lusztig \cite{LuGenerations} defined characters
\[ 
E_\l^0, E_\l^1, E_\l^2, \dots, E_\l^\infty \in (\ZM X)^\Wf
\]
for fixed $p$ and any highest weight $\l$. They are approximations to the character
of $L_\l$ in the following sense: $E_\l^0$ is the character of the simple
highest weight module in characteristic $0$; $E_\l^1$ is the
character of the simple highest weight module for a
quantum group at a $p^{th}$-root of unity; $E_\l^n$ is obtained from
$E_\l^{n-1}$ by a formula involving Kazhdan-Lusztig polynomials; and for fixed $\l$ as $n \to \infty$ 
the $E_\l^n$ stabilise to the character
  $E_\l^\infty$ predicted by \eqref{LCF} and Steinberg's tensor
  product theorem. Thus $E^\infty_\l$ gives the character of
  $L_\lambda$ for large $p$. One might hope that $E_\l^n$ is the
  character of a simple highest weight module for a quantum group like
  object which has an ``$n$-step Steinberg tensor product
  theorem''. For $\mathfrak{sl}_2$ such an object (for any $n$) has
  recently been proposed by Angiono \cite{Angiono}.
\end{remark}

\subsection{Lusztig's character formula and weight
  multiplicity} After stating his conjecture, Lusztig noticed that it
implied an interesting property of certain Kazhdan-Lusztig polynomials
attached to the (extended) affine Weyl group. Namely, their value
at 1 gives the dimension of a weight space in a simple finite
dimensional representation of the Langlands dual group (see
\eqref{eq:qwm} below).


Lusztig's idea was to interpret what his character formula
predicts for Frobenius twists of simple modules. Let us first
introduce some notation. Given $\mu \in X_+$ set
\[
I(\mu) = \{ s \in S \; | \; s(\mu) = \mu \}
\]
and let $\Wf^\mu$ denote the set of minimal coset representatives for
$\Wf/\Wf_{I(\mu)}$. Given $\mu \in X_+$ we define
\[
\sigma^\mu = \sum_{x \in \Wf^\mu} e^{x\mu}.
\]
We will need the identity
\begin{equation}
  \label{eq:WCI}
  \sigma^\mu = \sum_{x \in \Wf^\mu} \e_x\chi_{ \mu - \rho + x\rho}
\end{equation}
which follows (after a little thought) from Weyl's character formula.

Now fix $\l \in X_+$ and suppose that $p$ is large enough so that $L_\l = \Delta_\l$
is simple and $p\l$ satisfies Jantzen's condition. Thus the character
\[
\ch L_\l = \sum_{\mu \in X_+} (\dim \Delta_\l(\mu)) \sigma^\mu
\]
is given by Weyl's character formula. We also have
\[
L_{p\l} = \Delta_\l^\Fr
\]
and hence
\begin{align*}
  \ch L_{p\l} &= \sum_{\mu \in X_+}   (\dim \Delta_\l(\mu)) \sigma^{p\mu}
  \\
& \stackrel{\eqref{eq:WCI}}{=} \sum_{\mu \in X_+} (\dim \Delta_\l(\mu)) \sum_{x \in \Wf^\mu}
\e_x\chi_{p\mu - \rho + x\rho}
\end{align*}
We can rewrite this in terms of the $p$-dilated dot action of $\Wae$ as
\begin{align*}
 \ch L_{t_\l w_0 \pdot (-2\rho)} & = \sum_{\mu \in X_+} (\dim \Delta_\l(\mu)) \sum_{x \in \Wf^\mu}
\e_x \ch \Delta_{ t_\mu  x w_0 \pdot (-2\rho)}.
\end{align*}
Comparing this with \eqref{LCF} we arrive at the
prediction:\footnote{We use the version of \eqref{LCF} with the
extended affine Weyl group, see Remark \ref{rem:LCF}\eqref{rem:LCFextended}.}
\begin{equation}
  \label{eq:qwm}
  h_{t_\mu w_0, t_\l w_0}(1) = \dim \Delta_\l(\mu).
\end{equation}
In other words, \eqref{LCF} predicts that the dimensions of weight
spaces of Weyl modules occur as values at 1 of Kazhdan-Lusztig
polynomials for $\Wae$. This fact was proven by Lusztig \cite{Lq}
(independently of his conjecture).

\begin{remark}
Because Kazhdan-Lusztig polynomials
have non-negative coefficients this gives a refinement of the dimension of
the weight space. A representation theoretic interpretation of this
refinement was given by Brylinski \cite{Bryl} in terms of what is
nowadays called the Brylinski-Kostant filtration of the weight space.
\end{remark}


\subsection{Status of Lusztig's character formula} 

We give a brief summary of what is known about Lusztig's conjecture:
\begin{enumerate}
\item Lusztig formulated his conjecture in analogy to the 
  Kazhdan-Lusztig conjecture \cite{KLp}. At the time the case
  of $\SL_2$ was known and Jantzen had
  determined the characters of all simple modules for $\SL_3, \Sp_4,
  \Gtwo$ and $\SL_4$ using his sum formula
  \cite{JantzenForms}. Jantzen
  had also noticed that a character formula for large $p$ would also
  determine the characters of simple highest weight modules in
  characteristic zero \cite[Corollar im Anhang]{JantzenHWBook}. The fact that his conjecture implied
  the Kazhdan-Lusztig conjecture was also noticed by Lusztig.
\item The first proof of the Lusztig conjecture for $p \gg 0$ was obtained by
  combining works of Kashiwara-Tanisaki \cite{KT1,KT2} (relating Kazhdan-Lusztig
  polynomials and affine Lie algebras), Kazhdan-Lusztig
  \cite{KLaffine, KL2, KL3} (relating
  affine Lie algebras and quantum groups at a root of unity), Lusztig
  \cite{LuM,LuMe}
  (handling the non-simply-laced case) and
  Andersen-Jantzen-Soergel \cite{AJS} (relating quantum groups and modular
  representations of the Lie algebra). (These steps followed a program
  outlined by Lusztig in \cite{LusztigFD, LOQG}.)
 The main result of
  Andersen-Jantzen-Soergel is the existence of a finitely
  generated $\ZM$-algebra $B$ whose base change to a field of
  characteristic $p > h$ controls the
  principal block of restricted Lie algebra representations and whose base change
  to $\CM$ controls representations of the small quantum group at a
  root of unity. It is then possible to deduce Lusztig's conjecture
  for algebraic groups in characteristic $p \gg 0$ from the case of
  the quantum group. The algebra $B$ is not explicit and this
  approach did not yield any reasonable bounds on $p$. Over a decade later,
  a more direct route between perverse sheaves and the quantum group
  was provided by Arkhipov-Bezrukavnikov-Ginzburg \cite{ABG}.
\item More recently, Fiebig used his theory of moment graphs to
  provide a new proof of Lusztig's conjecture \cite{F}. The idea is to give a
  functor from a combinatorial category of ``moment graph sheaves'' associated to the affine
  Grassmannian to a combinatorial category constructed by
  Andersen-Jantzen-Soergel controlling Lie algebra
  representations. Lusztig's conjecture is then deduced from the
  decomposition theorem applied to intersection cohomology complexes
  on the affine Grassmannian. In essence, Fiebig's approach simplifies
  the original proof discussed above by providing a direct link between intersection
  cohomology complexes and the work of Andersen-Jantzen-Soergel. By a careful analysis of the combinatorics of moment
  graph sheaves Fiebig was able to give an explicit (enormous) bound above which Lusztig's
  conjecture holds \cite{F2} and establish the multiplicity one case
  \cite{F3}. Using recent work of Elias and the author establishing
  Soergel's conjecture \cite{EW2} and its local version \cite{WHL},
  all of the arguments used by Fiebig can be made entirely
  algebraic. Using related ideas, an algorithm using Soergel bimodules to produce the ``bad
  primes'' for Lusztig's conjecture was discovered by Libedinsky \cite{LLC}.
\item The localisation theorem \cite{BMR, BMR2} of Bezrukavnikov-Mirkovi\'c-Rumynin,
  provides a completely different approach to Lusztig's conjecture,
  which is closer to the original ($D$-module) proof of the
  Kazhdan-Lusztig conjecture. Working in the broader setting of Lie
  algebra representations the authors establish an equivalence of
  derived categories with coherent sheaves on
  Springer fibres. Roughly speaking, these categories are modules over
  the affine Hecke category \cite{BezRiche} and one can use an
  alternative realisation of this category \cite{Baffine} to deduce Lusztig's
  conjectures for Lie algebra representations for large
  $p$ \cite{BM}. These are known to imply Lusztig's character formula for $G_\bk$
  \cite{FiebigLCMG}.
\item Recently the author (building on joint work with Elias
  \cite{EW} and
  He \cite{HeW}) discovered many counter-examples to the expected bounds in
  Lusztig's conjecture \cite{WT,WIH}.
The upshot is that the above bounds (like $p
  \ge h$ or $p \ge 2h -3$) are much too optimistic. In fact, recent
  advances in number theory imply that there is no polynomial bound 
  in the Coxeter number for the validity of Lusztig's conjecture (see
  the appendix to \cite{WT} by Kontorovich, McNamara and the author). We will discuss
  these results in more detail in the next section.
\end{enumerate}

\section{Constructible Sheaves}\label{sec:sheaves}

\subsection{Notation} \label{sec:sheafb}
Let $X$ denote a complex algebraic variety acted
on by an algebraic group $H$. For simplicity we assume:
\begin{enumerate}
\item $H$ has finitely many orbits on $X$;
\item each orbit is simply connected;
\item each orbit is equivariantly simply connected; i.e. the
  stabiliser of any point in $H$ is connected.
\end{enumerate}
More generally, we will also allow $X$ to be an ind-variety with a
compatible action of a pro-algebraic group $H$, such that each finite
dimensional approximation satisfies the above conditions (see \cite[\S
2.7]{JMW2}).

For a fixed ring of coefficients $\bk$ we consider:
\begin{gather*}
D_H^b(X; \bk): \begin{array}{c} \text{the equivariant
    derived category} \\ \text{of constructible sheaves of
    $\bk$-modules on $X$,} \end{array} \\
D_{(H)}^b(X; \bk): \begin{array}{c} \text{the 
    derived category of sheaves} \\ \text{ 
constructible with respect to the $H$-orbits on $X$.} \end{array}
\end{gather*}
(In the ind-variety case these categories are defined as the direct
limits under extension by zero of the corresponding finite dimensional
approximations. In particular, any object has
finite dimensional support.)
We will sometimes ignore the coefficients and instead write $D_H^b(X)$ and
$D^b_{(H)}(X)$ if the context is clear.
We denote by $\For : D_H^b(X) \to D^b_{(H)}(X)$ the functor of
forgetting the equivariance. The full subcategories of perverse
sheaves are denoted
\begin{gather*}
  \Perv_{H}(X) \subset D^b_{H}(X), \\
  \Perv_{(H)}(X) \subset D^b_{(H)}(X).
\end{gather*}

Consider the decomposition of $X$ into $H$-orbits:
\[
X = \bigsqcup_{\l \in \Lambda} X_\l.
\]
Given $\l \in \Lambda$ we denote by $j_\l : X_\l \into X$ the
inclusion, $d_\l$ the (complex) dimension of $X_\l$ and by
$\underline{\bk}_{X_\l}$ the constant sheaf on $X_\l$. We have the perverse sheaves
\begin{gather*}
\Delta^\bk_\l := {}^p j_{\l !} (\underline{\bk}_{X_\l}[d_\l]), \quad
\ic^\bk_\l := {}^p j_{\l !*} (\underline{\bk}_{X_\l}[d_\l]) \quad
\nabla^\bk_\l := {}^p j_{\l *} (\underline{\bk}_{X_\l}[d_\l])
\end{gather*}
which by abuse of notation we regard as objects both of
$\Perv_{H}(X)$ and $\Perv_{(H)}(X)$ (it will be clear from the context
which object we mean below). As above, we will sometimes drop the superscript
indicating the coefficients if it is clear from the context. Note that if $j_\l$ is an affine
morphism (in particular if $X_\l$ is affine) then
\[
\Delta^\bk_\l = j_{\l !} (\underline{\bk}_{X_\l}[d_\l]), \quad
\nabla^\bk_\l = j_{\l *} (\underline{\bk}_{X_\l}[d_\l]). \]
Our assumptions on $X$ and $H$ guarantee that if $\bk$ is a field then
the set
\[
\{ \ic_x^\bk \; | \; x \in \Lambda \}
\]
is a complete set of representatives for the simple objects in
$\Perv_H(X; \bk)$ (resp. $\Perv_{(H)}(X;\bk)$).

Sometimes it will be useful to take integral coefficients
below. Consider the functor of extension of scalars
\begin{gather*}
 (-) \otimes^L_{\ZM} \bk : D_{(H)}(X;\ZM) \to D_{(H)}(X;\bk).
\end{gather*}
For any $\l$ as above the object $\ic_\l^\ZM \otimes_\ZM^L \bk$ is
perverse but not simple in general. However if we fix $\l$ and allow
$\bk$ to vary then $\ic_\l^\ZM \otimes_\ZM^L \bk$  will fail to be
simple in only finitely many (positive)
characteristics. For background on decomposition numbers for perverse
sheaves the reader is referred to \cite{decperv}.

In addition let us assume
\[
H^j(X_\l;\ZM) = H^j_H(X_\l;\ZM) = 0 \quad \text{for $j$ odd, and all $\l \in \Lambda$.}
\]
For any field $\bk$ and $\l \in \Lambda$ we denote by
\[
\EC_\l^\bk \in D^b_H(X) \quad \text{(resp. $ \in D^b_{(H)}(X)$)}
\]
the indecomposable parity sheaf (if it exists). Recall that $\EC_\l^\bk$ is
characterised up to isomorphism by the properties: it is supported on
$\overline{X_\l}$; its restriction to $X_\l$ is isomorphic to
$\underline{\bk}_{X_\l}[d_\l]$; and, its stalks and costalks vanish in
degrees of parity differing from that of $d_\l$. We denote by
\[
\Parity_H(X) \subset D^b_H(X) \quad \text{and} \quad \Parity_{(H)}(X)  \subset D^b_{(H)}(X)
\]
the full subcategories of parity complexes.

\subsection{The affine Grassmannian and flag variety} \label{sec:affineGrassFlag}

Let $G^\vee$ denote a complex reductive group whose root datum is dual
to the root datum of $G_{\bk}$. Let $T^\vee \subset G^\vee$ denote the
dual torus. Because 
$G_{\bk}$ is assumed simply connected, $G^\vee$ is of adjoint type.

To $G^\vee$ we associate its algebraic loop
group $G^\vee((t))$. Let $K = G^\vee[[t]]$ denote a maximal compact subgroup. We
consider the evaluation at $t = 0$ map
\[
\ev : K = G^\vee[[t]] \stackrel{t=0}{\longrightarrow} G^\vee
\]
and consider the Iwahori subgroup
\[
\Iw := \ev^{-1}(B^\vee)
\]
where $B^\vee$ is the Borel subgroup in $G^\vee$ containing $T^\vee$
and whose Lie algebra contains all characters in $\Phi_+^\vee$.

The \emph{affine flag variety} and \emph{affine Grassmannian} are the spaces
\begin{gather*}
  \Fl := G^\vee((t))/\Iw, \\
  \Gr := G^\vee((t))/K.
\end{gather*}
Both are $\CM$-schemes (of infinite type) in a natural way. In fact,
they are both ind-projective varieties. The natural projection
\begin{gather}
  \label{eq:FlGr}
  p : \Fl \to \Gr
\end{gather}
realises $\Fl$ as a $G^\vee/B^\vee$-bundle over $\Gr$.

We have
\begin{gather*}
  \pi_0(\Fl) = \pi_0(\Gr) = \Omega = X_*/\ZM \Phi^\vee
\end{gather*}
where $\Omega$ is the set of length-zero elements introduced in
\S\ref{sec:Wa}.

\begin{remark}
All connected components of $\Fl$ and $\Gr$ are isomorphic as
ind-varieties (even as $\Iw$-varieties). Each connected component is
isomorphic to the Kac-Moody flag variety associated to the extended
Cartan matrix of $G^\vee$.  
\end{remark}

The $\Iw$-orbits on $\Fl$ yields the Bruhat decomposition
\begin{gather}
  \label{eq:2}
  \Fl = \bigsqcup_{x \in \Wae} \Flx \quad \text{where} \quad \Flx :=
  \Iw \cdot x \Iw/\Iw
\end{gather}
and each $\Iw$-orbit is an affine space of dimension
\begin{gather}
  \label{eq:Fldim}
  \dim \Flx = \ell(x).
\end{gather}

Similarly, the $\Iw$-orbits on $\Gr$ give the Bruhat decomposition
\begin{gather}
  \label{eq:bruhat}
    \Gr = \bigsqcup_{x \in \Wae/W} \Grx \quad \text{where} \quad \Grx :=
  \Iw \cdot x \Iw/K
\end{gather}
and each cell in the Bruhat decomposition is isomorphic to an affine space
of dimension
\begin{gather}
  \label{eq:GrIwDim}
\dim \Grx = \ell(x_-) = \ell(x_+) - \ell(w_0)  
\end{gather}
where $x_-$ (resp. $x_+$) is the minimal (resp. maximal) element in the coset $xW$.

Recall that any element of $X$ is a coweight of $G^\vee$ which we can
regard as a point $t^\lambda \in G^\vee[t,t^{-1}] \subset G^\vee((t))$. These points
feature in the decomposition of $\Gr$ into $K$-orbits:
\begin{gather}
  \label{eq:iwasawa}
    \Gr = \bigsqcup_{\l \in X_+ } \Grl  \quad \text{where} \quad \Grl :=K \cdot t^\l K/K.
\end{gather}
Each cell is of dimension
\begin{gather}
  \label{eq:GrIwDim}
\dim K \cdot t^\l K/K  = 2 \langle \l, \rho \rangle.
\end{gather}
Moreover, each cell is an affine space bundle over a partial flag
variety and in particular is simply connected.

Let $\bk$ denote a fixed field of coefficients. With a slight
variation on the notation of the previous section denote by
 \[  \ic^\bk_{x, \Iw} \in D^b_\Iw(\Fl), \; \ic^\bk_{x, K} \in
  D^b_\Iw(\Gr) \text{ and } \ic^\bk_{\l} \in D^b_K(\Gr)
\]
the simple perverse sheaves.

Recall the following classical theorem of Kazhdan and Lusztig:

\begin{thm}[Kazhdan-Lusztig \cite{KLSchubert}] \label{thm:klpoly}
For $x, y\in \Wae$ we have
\[
h_{y,x} = \sum_{i \in \ZM} \dim H^{-i}( (\ic^\QM_{x,\Iw})_y) v^{
  i - \ell(y)}.
\]
\end{thm}

Because $h_{y,x} \in v^{\ell(x) -\ell(y)} \ZM[v^2]$, setting $v = -1$
we deduce that
\[
\e_x\e_y (h_{y,x}(1)) = \e_y \chi(( \ic_{x,\Iw}^{\QM})_y)
\]
or in other words
\begin{equation}
  \label{eq:ICchi}
  \chi(( \ic_x^\QM)_y) = \e_x h_{y,x}(1).
\end{equation}

We now deduce a similar formula for the affine
Grassmannian. Fix $x, y \in \Wae$ and assume that $x \in xW$ and $y
\in yW$ are maximal. Because $p : \Fl \to \Gr$ is a
smooth fibration with fibre $G^\vee/B^\vee$ and $\Flx\subset p^{-1} \Grx$ is open and dense
we have
\[
p^*\ic_{x, K}[\ell(w_0)] \cong \ic_{x,\Iw}.
\]
Thus
\begin{equation}
  \label{eq:ICKchi}
  \chi((\ic_{x,K}^\QM)_y) = \e_{w_0} \chi((\ic_{x,\Iw}^\QM)_{y} =
  \e_{w_0}\e_x h_{y,x}(1).
\end{equation}

\subsection{The (extended) affine Hecke category} \label{sec:heckecat}

Recall the affine Hecke algebra $\He$ and the extended affine Hecke algebra $\Hee$ 
from \S\ref{sec:hecke}. In this section we will describe
categorifications of $\He$ and $\Hee$ via sheaves on $\Fl$.

The category $D^b_{\Iw}(\Fl)$ is equipped with the structure of a
monoidal category:
\begin{gather*}
  * : D^b_{\Iw}(\Fl) \times D^b_{\Iw}(\Fl) \to D^b_{\Iw}(\Fl)
\end{gather*}
defined via
\begin{gather*}
  \FS * \GS := m_* ( \widetilde{\FS\GS})
\end{gather*}
where $m : G^\vee((t)) \times_\Iw \Fl \to \Fl$ is induced by the
multiplication map on $G^\vee((t))$ and, in a shorthand notation,
\[
 \widetilde{\FS\GS} := \res_{\Iw^4}^{\Iw^3} ( \FS \boxtimes \GS).
\]
In more detail, to construct $\widetilde{\FS\GS}$, we:
\begin{enumerate}
\item use the quotient equivalence $D^b_{\Iw \times \Iw} ( G^\vee((t)))
\simto D^b_\Iw( \Fl)$ to view the exterior tensor product $\FS
\boxtimes \GS \in D^b_{\Iw^4}( G^\vee((t)) \times G^\vee((t)))$;
\item restrict along the map $\Iw^3 \into \Iw^4:
(a, b, c) \mapsto (a, b^{-1}, b, c)$ to obtain an object 
$\res_{\Iw^4}^{\Iw^3}  (\FS
\boxtimes \GS) \in D^b_{\Iw^3}( G^\vee((t)) \times G^\vee((t)))$;
\item use the 
equivalence $D^b_{\Iw^3}( G^\vee((t)) \times G^\vee((t))) \simto
D^b_{\Iw}( G^\vee((t)) \times_{\Iw} \Fl)$ to view $\res_{\Iw^4}^{\Iw^3}(\FS \boxtimes \GS)$
as an object $\widetilde{\FS\GS} \in D^b_{\Iw}( G^\vee((t))
\times_{\Iw} \Fl)$.
\end{enumerate}

\begin{remark}
The above definition of
  convolution is mimicking convolution of
  $H$-biinvariant functions on a finite group $G$ (for $H \subset G$ a
  subgroup). It was in the context of Grothendieck's function-sheaf
  dictionary that this definition was first made \cite{Springer}.
\end{remark}

\begin{remark}
  The above definition makes sense if we regard $G^\vee((t))$ and
  $\Iw$ as (infinite-dimensional) topological groups. However if one
  wishes to work in a more algebraic category (for example to apply
  the decomposition theorem) then more care is needed, see \cite[\S
  2.2 and \S 3.3]{Nad}.
\end{remark}

\begin{remark} \label{rem:parity}
  Below an important property of convolution is that it preserves parity
  complexes \cite[Theorem 4.8]{JMW2}. In its basic form, this
  observation goes back to Springer, Brylinski and MacPherson
  \cite{Springer}. Its importance for modular representation theory was emphasised by Soergel in
  \cite{Soe}.
\end{remark}

Given a collection of objects $\{ \FC_i \}_{i \in I}$ in
$D^b_\Iw(\Fl)$ let $\langle \FC_i \rangle_{*, \oplus, [\ZM]}$ denote
the additive, graded, monoidal envelope of $\{ \FC_i \}_{i \in I}$:
its objects are the direct sums of shifts of tensor products
\[
\FC_{\textbf i} := \FC_i * \FC_j * \dots * \FC_k
\]
for any finite sequence ${\textbf i} = (i, j, \dots, k)$ of elements
of $I$. (We allow the empty sequence and set $\FC_{\emptyset} := \EC_{\id}$.)
We denote by $\langle \FC_i \rangle_{*, \oplus, [\ZM], \Kar}$ the category
obtained from $\langle \FC_i \rangle_{*, \oplus, [\ZM]}$ by adjoining all direct
summands of objects.

\begin{remark}
  The notation is intended to remind us that there is a formal
  procedure (``Karoubi envelope'') which allows us to produce (a
  category equivalent to) $\langle \FC \rangle_{*, \oplus, [\ZM],
    \Kar}$ starting from $\langle \FC_i \rangle_{*, \oplus,
    [\ZM]}$.
\end{remark}

We set
\begin{gather*}
  \HC_\bs := \langle \ES_s \; | \; s \in S \rangle_{*, \oplus, [\ZM]},
  \quad
  \HCe_\bs := \langle \{ \ES_s \; | \; s \in S\} \cup \{ \ES_\omega \;
  | \; \omega \in \Omega \} \rangle_{*, \oplus,
    [\ZM]},
\\
\HC := \langle \ES_s \; | \; s \in S \rangle_{*, \oplus, [\ZM], \Kar},
\quad
  \HCe := \langle \ES_s \; | \; s \in S\} \cup \{ \ES_\omega \;
  | \; \omega \in \Omega \} \rangle_{*, \oplus, [\ZM],
  \Kar}.
\end{gather*}
We call $\HC$ and $\HCe$ the \emph{(extended) affine Hecke
  category}. One may check that $\HC$ is the full subcategory of
$\HCe$ consisting
of objects supported on the identity component of $\Fl$.

It follows from Remark \ref{rem:parity} that every object of $\HCe$ is
parity. In fact, one has equalities
\begin{gather*}
  \HCe = \Parity_\Iw(\Fl) \quad \text{and} \quad
  \HCe = \Parity_\Iw((\Fl)_0)
\end{gather*}
where $(\Fl)_0$ denotes the identity component of $\Fl$. The following
theorem explains the name of $\HC$ and $\HCe$.

\begin{thm}
  There exists a unique isomorphism
\[
\Hee \simto [\HCe]_\oplus
\]
such that $\un{h}_s \mapsto [\ES_s]$ for all $s \in \Sa$ and $\omega
\mapsto \ES_\omega$ for all $\omega \in \Omega$. This isomorphism
induces an isomorphism $\He \simto [\HC]_\oplus$.
\end{thm}

 (Recall that $[\HCe]_\oplus$ denotes the split Grothendieck group of $\HCe$, see
\S\ref{sec:not}.) We denote by
\[
\ch : [\HCe_\oplus] \to \Hee
\]
the inverse to the
isomorphism of the theorem. It may be described explicitly via:
\[
\ch(\FS) = \sum_{x \in \Wae} \left ( \sum_{i \in \ZM} \dim H^{-i}( \FS_y) v^{
  i - \ell(y)} \right ) h_x.
\]

We define elements $\p\un{h}_x$ via
\[
\p\un{h}_x := \ch(\EC_x).
\]
Then $\{ \p\un{h}_x \}$ is a basis for $\Hee$ called the \emph{$p$-canonical
  basis}. (It only depends on the characteristic of $\bk$.) If we write
\[
\p\un{h}_x := \sum_{y \le x} \p h_{y,x} h_y 
\]
then the polynomials $\p h_{h,x}$ are the \emph{$p$-Kazhdan-Lusztig
  polynomials}. Here are some basic properties:
\begin{enumerate}
\item For fixed $x \in \Wae$ and $p \gg 0$ we have $\p\un{h}_x = \un{h}_x$.
\item If we write
\[
\p\un{h}_x = \sum \p a_{y,x} \cdot \p \un{h}_y
\]
then $\p a_{y,x} \in \ZM_{\ge 0}[v,v^{-1}]$ and $\overline{\p a_{y,x}}
= \p a_{y,x}$. In particular, $\p h_{y,x}$ have $\ge 0$ coefficients.
\item If we write
\[
\p\un{h}_x\p\un{h}_y = \sum \p \mu^z_{y,x} \cdot \p \un{h}_z
\]
then $\p \mu^z_{y,x} \in \ZM_{\ge 0}[v,v^{-1}]$ and $\overline{\p \mu^z_{y,x}}
= \p \mu^z_{y,x}$.
\end{enumerate}
For further discussion and examples the reader is referred to
\cite{JensenW}.

\begin{remark}
  Kazhdan-Lusztig polynomials depend only on the underlying Coxeter
  system. This is no longer true for $p$-Kazhdan-Lusztig polynomials
  (see \cite{JensenW}).
\end{remark}

\begin{remark}
  It is important for computations that $\HC$ and $\HCe$ have an alternative
  diagrammatic / algebraic presentation via generators and relations
  \cite{EW,RW}.
\end{remark}

\subsection{Geometric Satake equivalence} \label{sec:Satake}

As for $D^b_\Iw(\Fl)$, we can equip $D^b_K(\Gr)$ with the structure of
a monoidal category via
\[
\FS * \GS := m_* \widetilde{\FS\GS}
\]
where $m : G^\vee((t)) \times_K \Gr \to \Gr$ is induced by the
multiplication on $G^\vee((t))$ and
\[
 \widetilde{\FS\GS} := \res_{K^4}^{K^3} ( \FS \boxtimes \GS).
\]
is defined by mimicking the construction in \S \ref{sec:heckecat},
with $K$ in place of $\Iw$.

In this setting two miracles occur:
\begin{enumerate}
\item The convolution preserves $\Perv_{K}(\Gr)$:  if $\FS, \GS
  \in \Perv_K(\Gr)$ then so is $\FS * \GS$;
\item The convolution $*$ is symmetric: we have a canonical
  isomorphism $\FS * \GS \simto \GS * \FS$ equipping $\Perv_K(\Gr)$ with
  the structure of symmetric tensor category.
\end{enumerate}

Recall our semi-simple algebraic group $G_{\bk}$ over $\bk$ whose root
system is dual to that of $G^\vee$.

\begin{thm}[Geometric Satake \cite{MV}] There is an equivalence of monoidal 
  categories
\[
  \SC : ( \Perv_K(\Gr; \bk), *) \simto (\Rep G_{\bk}, \otimes_{\bk} ).
\]
\end{thm}

More generally, this theorem is true with coefficients in any
commutative ring.\footnote{In \cite{MV} the result is proved under the
  additional assumption that the ring of coefficients if Noetherian
  and of finite global dimension. However this is for simplicity only
  (see the discussion at the bottom of pg. 100 of \cite{MV}).}

\begin{remark} Some remarks on the geometric Satake equivalence:
  \begin{enumerate}
 \item A remarkable aspect of the proof of geometric Satake is that
    it does not construct a functor in either direction! Instead, one
    uses the Tannakian formalism to deduce that $\Perv_K(\Gr; \ZM)$ is
    equivalent to the representations of some group scheme 
    over $\ZM$. After considerable work, one manages to identify this
    group scheme with the Chevalley group scheme. Thus the equivalence
    can actually be seen as providing a
    \emph{construction} of the dual group. In the above notation, $G_\ZM$
    is
    constructed starting from $G^\vee$. As far as the author is
    aware, geometric Satake is the only known construction of the
    Langlands dual group.
  \item There is no proof of geometric Satake which works
    directly with coefficients in characteristic $p$. At present, the
    case of coefficients of characteristic $p$ is deduced by reduction
    modulo $p$ from the corresponding statement over $\ZM$.
  \item Recent work of Mautner-Riche \cite{MR} would yield a new proof
    that $\SC$ is an equivalence of abelian categories if one could prove that $\Perv_K(\Gr;k)$
    is a highest weight category without using geometric Satake.
  \item Most of the basic theorems concerning the representation
    theory of $G_\bk$ have no proof on the side of perverse sheaves. For
    example, at present there is no geometric proof of either the
    Steinberg tensor product theorem or the linkage principle.
  \item If the characteristic $p$ of $\bk$ is good then the parity
    sheaves $\EC_\l \in D^b_K(\Gr; \bk)$ are perverse and correspond
    under the geometric Satake equivalence to tilting modules
    \cite{JMW3, MR}. This gives a proof via constructible sheaves
    of Theorem \ref{thm:ttt} (at least in good characteristic). This
    result was used by Achar and Rider to show that the stalks of
    $\Delta_\l^\ZM$ are free of $p$-torsion if $p$ is good for $G$
    (``Mirkovi\'c-Vilonen conjecture'' \cite{Ju-AffGr}) \cite{ARMV}.
  \end{enumerate}
\end{remark}

\subsection{The Finkelberg-Mirkovi\'c conjecture} \label{sec:FM}
The geometric Satake
equivalence is not so useful for studying character formulas, because
the linkage principle is not visible. (Indeed, as we have tried to
explain, it is more often the case that theorems on representations of $G_{\bk}$
predict remarkable behaviour on the constructible side of
geometric Satake which have no geometric explanation at present.)
In this section we explain a conjecture of Finkelberg-Mirkovi\'c which
gives a geometric realisation of the extended principal block. Because it already
incorporates the linkage principle, it is much more useful for
understanding characters (amongst other things).

Recall that $\Rep$ denotes the category of finite-dimensional algebraic representations
of our semi-simple group $G_{\bk}$, which is defined over a field
$\bk$. Also recall the extended principal block $\Repe \subset
\Rep$. In \S\ref{sec:repp} we explained that Frobenius twist makes
$\Repe$ a right module category over $\Rep$ via $(V, M)
\mapsto V \otimes M^\Fr$.

By imitating the convolution product on $\Perv_K(\Gr)$ one can define
a convolution product
\[
* : D^b_{(\Iw)}(\Gr) \times D^b_K(\Gr) \to D^b_K(\Gr)
\]
which again (miraculously) descends to a convolution product
\[
* : \Perv_{(\Iw)}(\Gr) \times \Perv_K(\Gr) \to \Perv_K(\Gr),
\]
making $\Perv_{(\Iw)}$ a right module category over
$\Perv_K(\Gr)$.

\begin{conj}[Finkelberg-Mirkovi\'c \cite{FM}]
  There is an equivalence
\[
\QC : \Repe \simto \Perv_{(\Iw)}(\Gr, k) 
\]
mapping
\begin{gather*}
  L_{x \pdot (-2\rho)} \mapsto \ic_{x^{-1}} \\
  \Delta_{x \pdot (-2\rho)} \mapsto \Delta_{x^{-1}}
\end{gather*}
for all $x \in \Wae$ such that $x \pdot (-2\rho) \in X_+$.

Moreover, this equivalence is compatible with geometric Satake and
Frobenius twist. That is, the following diagram is commutative
up to natural isomorphism:
\begin{equation}
  \label{eq:FMdiag}
\begin{array}{c}
\tikz[xscale=2,yscale=0.6]{
\node (ul) at (-1,1) {$\Repe \times \Rep$};
\node (ur) at (1,1) {$\Repe$};
\draw[->] (ul) to node[above]{\tiny $(-) \otimes (-)^{\Fr}$} (ur);
\node (ll) at (-1,-1) {$\Perv_{(\Iw)}(\Gr) \times \Perv_K(\Gr)$};
\node (lr) at (1,-1) {$\Perv_{(\Iw)}(\Gr)$};
\draw[->] (ll) to node[above]{\tiny $(-) * (-)$} (lr);
\draw[->] (ul) to node[right]{\tiny $\QC \times \SC$} (ll);
\draw[->] (ur) to node[right]{\tiny $\QC$} (lr);
}
\end{array}
%
\end{equation}
\end{conj}

\begin{remark} Some remarks concerning the Finkelberg-Mirkovi\'c conjecture:
  \begin{enumerate}
  \item At first sight the natural inclusions on both sides of this
    equivalence appear to go in opposite (``wrong'') directions: on the
    representation theory side we  have $\Repe \subset \Rep$; whereas on the
    perverse sheaf side forgetting
    $K$-equivariance defines a fully-faithful embedding $\For_K : \Perv_K(\Gr)
    \into \Perv_{(\Iw)}(\Gr)$. This ``contradiction'' is resolved by
    Frobenius twist. By acting on the right on the trivial
    representation, the commutativity of
    \eqref{eq:FMdiag} implies that we have a commutative diagram up
    to natural isomorphism:
\begin{equation}
  \label{eq:FMdiag2}
\begin{array}{c}
\tikz[xscale=2,yscale=0.6]{
\node (ul) at (-1,1) {$\Rep$};
\node (ur) at (1,1) {$\Repe$};
\draw[right hook-latex] (ul) to node[above]{\tiny $(-)^{\Fr}$} (ur);
\node (ll) at (-1,-1) {$\Perv_K(\Gr)$};
\node (lr) at (1,-1) {$\Perv_{(\Iw)}(\Gr)$};
\draw[right hook-latex] (ll) to node[above]{\tiny $\For_K$} (lr);
\draw[->] (ul) to node[right]{\tiny $\SC$} (ll);
\draw[->] (ur) to node[right]{\tiny $\QC$} (lr);
}
\end{array}
\end{equation}
In fact, \cite{FM} ask only for the commutativity of
\eqref{eq:FMdiag2}.
\item The Finkelberg-Mirkovi\'c conjecture is still a conjecture. Recently, Achar and Riche have come very
  close to a proof \cite{ARFM}. Building on work of Achar-Riche \cite{ARMixed1,ARMixed2,ARMixed3},
  Achar-Rider \cite{ARMV,AR} and Mautner-Riche \cite{MR} they prove that a
  certain ``mixed'' version of the category $\Perv_{(\Iw)}(\Gr)$
  provides a graded enhancement of $\Repe$, and check compatibility
  with geometric Satake. The remaining difficulty is to construct a
  ``forgetting the mixed structure'' functor to $\Perv_{(\Iw)}(\Gr)$.
\item The Finkelberg-Mirkovi\'c conjecture predicts
  that $\Perv_{(\Iw)}(\Gr, \ZM)$ provides an abelian category over
  $\ZM$ from
  which the principal block of $G_{\bk}$ in any characteristic $p \ge h$ may
  be deduced by ``reduction modulo $p$''. This provides a simple
  explanation for many independence of $p$ results, as we hope will
  become clear. The
  existence of such an integral form seems mysterious from an
  algebraic point of view.
  \end{enumerate}
\end{remark}

Let us explain how the Finkelberg-Mirkovi\'c conjecture may be used to
deduce character formulas for $\Repp$ in terms of the geometry of the
affine Grassmannian. As usual, we start by writing
\[
[L_{x \pdot (-2\rho)}] = \sum a_{y,x} [\Delta_{y \pdot (-2\rho)}]
\]
for certain $a_{y,x} \in \ZM$. Applying the equivalence $\QC$ from the
Finkelberg-Mirkovi\'c conjecture we deduce
\[
[\ic_{x^{-1}}] = \sum a_{y,x} [\Delta_{y^{-1}}].
\]
Taking the Euler characteristic on both sides at a point of the stratum
$\EuScript Gr_{y^{-1}}^\vee$ 
yields\footnote{In case the reader wants to worry about signs: The dimension of $\EuScript Gr_{y^{-1}}^\vee$  is $\ell(y^{-1}) -
\ell(w_0)$, because $y \in Wy$ is maximal. This is where the
$\e_{w_0}\e_{y^{-1}}$ comes from.}
\begin{equation}
  \label{eq:ayx}
a_{y,x} = \e_{w_0}\e_{y^{-1}}\chi( (\ic_{x^{-1}})_{y^{-1}}).  
\end{equation}
Thus the $a_{y,x}$ are (up to sign) simply the Euler characteristics of
the stalks of the intersection cohomology complexes on $\Gr$!

Moreover, if the characteristic $p$ of $\bk$ is large enough then,
\[
\chi( (\ic^\bk_{x^{-1}, K})_{y^{-1}}) = \chi(
(\ic^\QM_{x^{-1}, K})_{y^{-1}}) \stackrel{ \eqref{eq:ICKchi}}{=} \e_{x^{-1}}\e_{w_0}
h_{y^{-1},x^{-1}}(1)
\]
or in other words (using that $h_{y^{-1},x^{-1}} = h_{y,x}$ and
$\e_{x^{-1}} = \e_x, \e_{y^{-1}} = \e_y$)
\begin{equation}
  \label{eq:Eulermult}
a_{y,x} = \e_{yx} h_{y,x}(1) \quad \quad \text{for $p$ large}.  
\end{equation}
This is the prediction made by Lusztig's conjecture. Recall that in
order to confirm Lusztig's conjecture we only need to check
\eqref{eq:Eulermult} for finitely many $x$ and $y$ (see Remark
\ref{rem:affine_finite}). Thus the Finkelberg-Mirkovi\'c conjecture implies Lusztig's conjecture
for large $p$, and helps us have some picture about what might ``go wrong''.

\subsection{Lusztig's conjecture and torsion} In the previous
section we have explained why the Finkelberg-Mirkovi\'c conjecture
implies character formulas in terms of stalks of intersection
cohomology complexes with coefficients in $\bk$. In this section we
explain (still assuming the Finkelberg-Mirkovi\'c conjecture) that
torsion in local integral intersection cohomology controls Lusztig's
conjecture.

Fix $p \ge h$ and consider the subsets of $\Wa$ defined as follows:
\begin{gather*}
  J_p := \{ x \in \Wa \; | \; \langle \a^\vee, x^{-1} \pdot (-2\rho) +
  \rho \rangle \le
p(p-h+2) \text{ for all $\a \in \Phi_+$}\}, \\
R := \{ x \in \Wa \; | \; x^{-1} \pdot (-2\rho) \in X_1^p \}.
\end{gather*}
The set $R$ is independent of $p$.

\begin{thm} \label{thm:LCFtorsion}
  Assume $p \ge h$ and the Finkelberg-Mirkovi\'c conjecture.
  \begin{enumerate}
  \item The original form of Lusztig's conjecture (Conjecture
    \ref{LC:original}) holds if and only if
    $\ic_{x,K}^\ZM$ has no $p$-torsion in its stalks and costalks, for all $x$
    in $J_p$.
\item The revised form of Lusztig's conjecture (Conjecture
    \ref{LC:revised}) holds if and only if
    $\ic_{x,K}^\ZM$ has no $p$-torsion in its stalks and costalks, for all $x$
    in $R$.
  \end{enumerate}
\end{thm}

\begin{remark}
We have explained in \S\ref{sec:LCSteinberg} why Kato's theorem implies that the equivalent conditions in
(2) imply the equivalent conditions in (1) if $p \ge 2h-3$. From a
geometric perspective this is rather surprising: for large $p$
there are many more Schubert varieties in $\Gr$ parametrised by $J_p$ than by $R$.
\end{remark}

\begin{remark} \label{rem:ext}
Taking into account Remark \ref{rem:LCF}(1), identical arguments to those below show
  that Lusztig's conjecture is equivalent to the absence of
  $p$-torsion in the stalks or costalks of $\ic_{x,K}^\ZM$ for any $x$
  belonging to the set
\[
 J^\ext_p := \{ x \in \Wae \; | \; \langle \a^\vee, x^{-1} \pdot (-2\rho) +
  \rho \rangle \le
p(p-h+2) \text{ for all $\a \in \Phi_+$}\}
\]
\end{remark}

The following is an immediate consequence of (2) and the above remark:

\begin{cor} Let $\kappa \ge 2h-3$ and assume the Finkelberg-Mirkovi\'c conjecture. Suppose that all stalks and  costalks of
  $\ic_{x,K}^\ZM$ are free of $p$-torsion, for all $p \ge \kappa$ and all $x
  \in R$. Then
  (both formulations of) Lusztig's conjecture hold in all characteristics $p \ge \kappa$.
\end{cor}

Let us explain why Theorem \ref{thm:LCFtorsion} holds. Suppose first
that the stalks and costalks of $\ic_{x,K}^{\ZM}$ are free of
$p$-torsion, for all $x \in J_p$ (resp. $R$). Then
\[
\ic_{x,K}^{\ZM} \otimes_{\ZM}^L \bk = \ic_{x, K}^\bk
\]
and hence
\[
\chi((\ic_{x,K}^{\QM})_y) = \chi((\ic_{x,K}^{\ZM} \otimes_{\ZM}^L \bk)_y) = \chi((\ic_{x, K}^\bk)_y)
\]
for all $x \in J_p$ (resp. $x \in R$) and all $y$. This implies the
original (for $x \in J_p$) and revised (for $x \in R$) forms of Lusztig's conjecture, as we have explained in the
previous section.

The other direction is a little more involved. Assume that the stalks and costalks of
$\ic_x^{\ZM}$ are not free of $p$-torsion, for some $x \in J_p$
(resp. $x \in X$). Hence $\pcan_x \ne \un{h}_x$ (see
\cite[Corollary 3.13]{W}). If we assume in addition that $x \in J_p$
(resp. $X$) is minimal such that $\pcan_x \ne \un{h}_x$, then ${}^p
a_{y,x} \in \ZM$ for all $y \le x$ (see \S\ref{sec:heckecat} for the
notation $a_{y,x}$). Now, by
\cite[Proposition 2.3]{WilRed} it follows that $\ic_{x,K}^{\ZM}
\otimes_\ZM^L \bk$ has a non-trivial decomposition number; in other
words that
\[
[\ic_{x,K}^{\ZM} \otimes_\ZM^L \bk] \ne [\ic_{x,K}^{\bk}] \quad
\text{in $[\Perv_{(\Iw)}(\Gr, \bk)]$.}
\]
Thus there exists a $y$ such that 
\[
\chi((\ic_{x,K}^{\QM})_y) = \chi((\ic_{x,K}^{\ZM} \otimes_{\ZM}^L \bk)_y) \ne \chi((\ic_{x, K}^\bk)_y)
\]
and \eqref{LCF} cannot hold for the simple module $L_{x^{-1} \pdot
  (-2\rho)}$, as we explained in the previous section.

\begin{remark}
  The implication $\Leftarrow$ of Theorem \ref{thm:LCFtorsion}(2) is
  a theorem of Fiebig \cite{F} (see also \cite{FW}). It is a key
  ingredient in his proof of Lusztig's conjecture for large $p$, as
  well as his bound \cite{F2}.
\end{remark}

\begin{remark}
One needs to be careful when dealing with the set $R$ because
its image in $\Wa / W$ is usually not closed in the Bruhat
order. (That is, there exist Bruhat cells $\Grx$ for $x \in
R$ whose closures contain cells $\EuScript Gr_{y}^\vee$ with $y
\notin R$.)
\end{remark}

\subsection{Counter-examples to expected bounds} 

In this section we describe the results of \cite{WT,WIH} in the context of
the Finkelberg-Mirkovi\'c conjecture.

In the previous section we explained (assuming the
Finkelberg-Mirkovi\'c conjecture) that torsion in the stalks and
costalks of $\ic_{x,K}^\ZM$ for certain $x \in \Wa$ controls
Lusztig's conjecture. However it seems very difficult to compute or understand the
torsion in the stalks of $\ic_{x,K}^\ZM$. Soergel suggested that it might
be worthwhile to study the finite flag variety as a ``toy model'' for
the geometric study of Lusztig's conjecture.\footnote{In the words of
  Soergel \cite{Soe}: ``The goal of this article is to forward this
  problem [Lusztig's conjecture] to the topologists or geometers.''}

To explain the relevance of the finite flag variety we need to recall
the notion of smooth equivalence. A \emph{singularity} is a pair $(X,x)$ where $X$ is an algebraic
variety and $x \in X$ is a point. Two singularities $(X,x)$ and
$(Y,y)$ are \emph{smoothly equivalent} if there exists another
singularity $(Z,z)$ and smooth maps $X \stackrel{f}{\leftarrow} Z
\stackrel{g}{\to} Y$ with $x = f(z)$
and $y = g(z)$. If $(X,x)$ and $(Y,x)$ are smoothly equivalent then
small Euclidean neighbourhoods of $x$ and $y$ are analytically
isomorphic, up to taking a product with a smooth
variety. In particular, the stalks $(\ic_X^\ZM)_x$ and $(\ic_Y^\ZM)_y$
are isomorphic up to shift (see e.g. \cite[Proposition 3.8]{decperv}).

Consider the finite flag variety $X^\vee := G^\vee/B^\vee$ and its
Bruhat decomposition
\[
X^\vee = \bigoplus_{x \in W} X^\vee_{x} \quad \text{where $X^\vee_x :=
  B^{\vee} \cdot xB^\vee/B^\vee.$}
\]
We denote by $\ic_{x,B^\vee}^\ZM$ the integral intersection complex of the Schubert variety
$\overline{X^\vee_x}$.

Suppose that $z \in \Wae$ is minimal in its coset $Wz$ and choose elements $x =
x'z$ and $y = y'z$ with $x', y' \in W$ (so $\ell(x) =
\ell(x') + \ell(z)$, $\ell(y) = \ell(y') + \ell(z)$). We have:\footnote{A sketch: the map $gB^\vee \mapsto g \cdot z\Iw /\Iw$ defines an embedding
$X^\vee \into \Fl$. Consider $E = \bigsqcup_{x \in Wz} \Flx$. There exists a morphism $E \to X^\vee$ making $E$ into an affine space
bundle over $X^\vee$. This map is compatible with the $\Iw$
(resp. $B^\vee$) orbits on $E$ and $X^\vee$ and induces a smooth
morphism $\overline{\Flx} \cap E \to \overline{X^\vee_{x'}}$.}
\begin{gather*}
  \text{The singularities $(\overline{\Flx}, y)$ and
    $(\overline{X^\vee_{x'}},y')$ are smoothly equivalent.}
\end{gather*}
Moreover, if $x$ and $y$ are maximal in their cosets $xW$ and $yW$
respectively then (by considering the smooth map $p : \Fl \to \Gr$) we
conclude:
\begin{gather*}
  \text{The singularities $(\overline{\Grx}, y)$, $(\overline{\Flx}, y)$ and
    $(\overline{X^\vee_{x'}},y')$ are smoothly equivalent.}
\end{gather*}
The upshot is that if $x$ and $y$ belong to the same right $W$-coset
then torsion in the stalk of $(\ic_{x,K}^\ZM)$ at $yK/K \in \Gr$ can
be calculated on the finite flag variety.

Consider the coset $Wt_{-\rho} w_0 \subset \Wae$. For $x \in W$ we have
\[
(xt_{-\rho} w_0)^{-1} \pdot (-2\rho) = w_0( p(-\rho) + x^{-1} \wdot
(-2\rho)) = p\rho + w_0x^{-1} \wdot (-2\rho).
\]
We conclude that $Wt_{-\rho}w_0 \subset   J_p^{\ext}$ (see Remark
\ref{rem:ext}) if and only if
\[
\langle p\rho + \rho, \a^\vee \rangle \le p(p-h+2)  \text{ for all $\a \in \Phi_+$}.
\]
One may check that this is the case if $p >
2h-3$. From the results of the previous section we conclude:

\begin{thm} \label{thm:tornotor}
  If Lusztig's character formula is true for all $p \ge \kappa > 2h-3$ then
  there is no $p$-torsion in the stalks or costalks of the integral
  intersection cohomology complexes $\ic_{x,B^\vee}^\ZM$ for all $x
  \in W$ and $p \ge \kappa$.
\end{thm}

\begin{remark}
  The discussion above deduced the Theorem \ref{thm:tornotor} above
  from the
  Finkelberg-Mirkovi\'c conjecture. This is ahistorical, and the above
  theorem is known independently of the Finkelberg-Mirkovi\'c conjecture. In \cite{Soe}
  Soergel proves a theorem very similar to the the above formulation. The exact
  formulation above may be deduced by combining Soergel's results with
  the theory of parity sheaves \cite{JMW2}.
\end{remark}

For $m \ge 1$ define $T(m)$ to be the maximal prime number $p$ which
occurs as torsion in the stalk or costalk of some integral
intersection cohomology complex on the flag variety of $\GL_m(\CM)$
(if there is no torsion we set $T(m) = 1$). Because Schubert varieties
for $\GL_m$ are also Schubert varieties for $\GL_{m+1}$ our function
$T$ is monotonically increasing. Here is a table of some known values
of our function
\[
\begin{tabular} {c|ccccccccccccc}
$m$ & 1 & 2 & 3 & 4 & 5 & 6 & 7 & 8 & 9 & 10&  11&  12 & \dots\\ \hline
$T(m)$ & 1 & 1 & 1 & 1 & 1 & 1 & 1 & 2 & 2 & $\ge 2$ & $\ge 2$ & $\ge 3$ &
\dots
\end{tabular}
\]
(The values of $T(m)$ for $m \ge 9$ are due to Braden and the author:
for $T(m)$ with $m \le 8$ see \cite{W}, the value $T(9) = 2$ is
unpublished. The value $T(12) \ge 3$ is due to Polo (unpublished) who showed more
generally that $T(4p) \ge p$ for all primes $p$.)

Given the above (admittedly rather limited) data the following is
surprising:

\begin{thm}The function $T(m)$ grows at least exponentially in $m$.
\end{thm}

Combining this with the above results one obtains:

\begin{cor} Suppose that $\kappa(h)$ is a function of the Coxeter number,
  such that Lusztig's character formula holds for any $G_\bk$ in
  characteristic $p \ge \kappa(h)$. Then $\kappa(h)$ grows at least
  exponential in $h$.
\end{cor}

The main idea of \cite{WT} is that certain structure constants occurring in
Schubert calculus for the cohomology ring $H^*(X^\vee;\ZM)$ also occur as torsion in local
integral intersection cohomology in much higher rank groups. Using
these ideas it is shown, for
example, that any prime number dividing any entry of a word of length
$\ell$ in the semi-group
\[
{\small \left \langle \left ( \begin{matrix} 1 & 1 \\ 0 & 1 \end{matrix} \right ), 
 \left ( \begin{matrix} 1 & 0 \\ 1 & 1 \end{matrix} \right ) \right
\rangle } \subset \SL_2(\ZM)
\]
occurs as torsion in the stalks or costalks of some
$\ic_{x,B^\vee}^\ZM$ on the flag variety of $\GL_{5 +
  3\ell}(\CM)$. Some non-trivial number theory (which relies on recent
advances in ``thin groups'') gives the above results on torsion
growth.

\begin{remark}
  The main result of \cite{WT} uses a formula for certain entries of
  intersection forms obtained by the author and He \cite{HeW}. This
  result uses the theory of generators and relations for Soergel
  bimodules \cite{EW} in a crucial way. A purely geometric proof \cite{WIH} of the
  main result of \cite{WT} was
  discovered later.
\end{remark}

\subsection{Tilting modules and the Hecke category} In this section we
give a brief description of the conjectures and results of
\cite{RW}. The goal is to describe $\Rep_0$ (or more precisely its
tilting modules) via the Hecke
category. For more detail on any of the material below, the reader is
referred to \cite{RW}.

Recall that $\f\Wae$ denotes the set of minimal coset representatives
for $W \setminus \Wae$. It will be convenient to simplify notation as follows:
\[
L_x := L_{x \pdot 0},  \quad
\Delta_x := \Delta_{x \pdot 0}, \quad
\nabla_x := \nabla_{x \pdot 0}, \quad
T_x := T_{x \pdot 0} \quad \in \Rep_0.
\]

 Throughout we assume that $p > h$, where $h$ is the Coxeter
number. This choice guarantees that for all $s \in \Sa$ we may fix a weight
$\mu_s \in C_-^p$ whose stabiliser under the $p$-dilated dot action is
precisely $\langle s \rangle \subset \Wa$ (see \cite[\S
  6.3(1)]{JaBook}). We define the
\emph{wall-crossing functor} associated to $s \in \Sa$ as
\[
\Theta_s := T_{\mu_s}^{-2\rho} \circ T_{-2\rho}^{\mu_s} : \Rep_0 \to
\Rep_0.
\]
It will be convenient to view wall-crossing functors as acting on the
right.

Consider the \emph{anti-spherical module}
\[
\AntiS := \sgn \otimes_{\ZM W} \ZM \Wa = \bigoplus_{x \in \f\Wa} \ZM
\e \otimes x
\]
obtained by inducing the sign representation $\sgn = \ZM \e$ of the finite Weyl group
to the affine Weyl group. Because the classes $[\Delta_x]$ for $x \in \f\Wa$ span the
Grothendieck group we have an isomorphism:
\begin{align*}
\AntiS & \simto [\Rep_0], \\
\e \otimes x & \mapsto [\Delta_x] \quad \text{for all $x \in \f\Wa$.}
\end{align*}
Moreover,  it is an easy consequence of \cite[Chapter 7]{JaBook} that
we can upgrade to an isomorphism of right $\ZM \Wa$-modules if we make $[\Rep_0]$ into a
$\Wa$-module via
\[
[M] \cdot (1 + s) := [M\Theta_s] \quad \text{for all $s \in \Sa$.}
\]
Thus the action of wall-crossing functors on the
principal block categorifies the anti-spherical module.

The main conjecture of \cite{RW} is that this action of the affine
Weyl group on the Grothendieck group can be lifted to the Hecke
category:

\begin{conj} \label{conj:rw}
  $\Rep_0$ is a right module category over $\HC$, with $\EC_s$ acting via
  $\Theta_s$.
\end{conj}

\begin{remark}
  Actually, this is a slight simplification of the conjecture, which
  nonetheless captures its spirit. (The
  version in \cite{RW} requires that certain
  generating morphisms in $\HC$ arise from adjunctions between
  translation functors, see \cite[\S 5.1]{RW}.)
\end{remark}

\begin{remark}
  In \cite{RW} the above conjecture is proved for $\GL_n$. The proof
  uses the Chuang-Khovanov-Lauda-Rouquier theory of categorification of Lie
  algebras \cite{CR,R2KM,KhLa1,KhLa2}. (This theory is only available
  at present in type $A$.) It also makes essential use of a recent
  theorem of Brundan \cite{Brundan}.
\end{remark}

The main point of \cite{RW} is that the above conjecture has strong
structural and numerical consequences for $\Rep_0$. To discuss this
we need to explain another categorification of the anti-spherical
module.

The anti-spherical module is quantized via the right $\He$-module
\[
\AntiS_v := \sgn_v \otimes_{\He_f} \He = \bigoplus_{x \in \f\Wa} \ZM[v^{\pm 1}] n_x
\]
where $\sgn_v$ denotes the sign representation of the finite Hecke
algebra $\He_f$ given by $h_s \mapsto -v$ for
all $s \in S$ and $n_x := 1 \otimes h_x$ for $ x \in \f\Wa$. The
module $\AntiS_v$ has a canonical basis $\{ \un{n}_x \; | \; x \in
\f\Wa \}$ constructed by Deodhar \cite{Deodhar} (see also \cite{SoeKL}). We have:
\begin{gather*}
  n_{\id} \cdot \un{h}_x := \begin{cases} \un{n}_x & \text{if $x \in
      \f\Wa$,} \\
0 & \text{otherwise.} \end{cases}
\end{gather*}

It is not difficult to see that $\AntiS_v$ has alternative descriptions as
\[
\AntiS_v = \He / \langle \un{h}_s \He \; | s \in S \rangle = \He / (
\bigoplus_{x \notin \f\Wa} \ZM[v^{\pm 1}] \un{h}_x).
\]
Thus it is natural to try to categorify the anti-spherical quotient as
a quotient of additive categories:
\[
\CAS := \HC / \langle \EC_x \; | \; x \not\in \f\Wa
\rangle_{\oplus, [\ZM]}.
\]
(That is, $\CAS$ is defined to be the quotient of $\HC$ by the
ideal of morphisms factoring through any direct sum of shifts of
$\EC_x$, for some $x \not\in \f\Wa$.)

It is not difficult to see that $\CAS$ is a right module category
over $\HC$ and that the identification $\He = [\HC]_\oplus$
induces a canonical identification
\[
\AntiS_v = [\CAS]_\oplus
\]
of right $\He$-modules. The image of $\EC_x$ in $\CAS$ is
indecomposable if $x \in \f\Wa$ and is zero otherwise. Its class
\[
\p\un{n}_x := [\EC_x] \in \AntiS_v
\]
defines the \emph{$p$-canonical basis} in the anti-spherical
module. The $p$-canonical basis gives rise to the \emph{anti-spherical $p$-Kazhdan-Lusztig polynomials} $\p n_{y,x}$ via
\[
\p\un{n}_x := \sum_{y \in \f \Wa} \p n_{y,x} n_y.
\]
The basis $\{ \p \un{n}_x\}_{x \in \f\Wa}$ enjoys positivity
properties analogous to those of the $p$-canonical basis (see 
\S\ref{sec:heckecat}).

Recall that $\CAS$ is an additive graded category. Let us denote
by $\CAS_{/[\ZM]}$ the category obtained by ``forgetting the
grading'': it has the same objects as $\CAS$ and morphisms are given
by
\[
\Hom_{\CAS_{/[\ZM]}}(\EC, \EC') := \bigoplus_{m \in \ZM}
\Hom_{\CAS}(\EC, \EC'[m]).
\]
Let $\Tilt_0 \subset \Rep_0$ denote the full subcategory of tilting
modules. Note that $\Tilt_0$ is preserved by wall-crossing
functors.\footnote{This follows because translation functors preserve the categories of modules with good or Weyl filtration. Alternatively, one may
  appeal to Theorem \ref{thm:ttt} and the fact that we may choose tilting
modules to define our translation functors.} Thus if Conjecture
\ref{conj:rw} holds then $\Tilt_0$ is preserved by the action of the
Hecke category.

\begin{thm}
  Assume Conjecture \ref{conj:rw} holds:
  \begin{enumerate}
  \item We have an equivalence
\[
\CAS_{/[\ZM]} \simto \Tilt_0
\]
of $\HC$-module categories.
\item For all $x, y \in \f\Wa$ we have:
\[
(T_x:\Delta_y) = \p n_{y,x}(1).
\]
  \end{enumerate}
\end{thm}

\begin{remark} Some remarks concerning the theorem (assuming
  Conjecture \ref{conj:rw}):
  \begin{enumerate}
  \item Part (1) of the theorem implies that $\Tilt_0$ admits a
    grading (given by $\CAS$). In \cite{RW} it is explained how this
    grading can be used to produce a grading on $\Rep_0$. Another
    grading on $\Rep_0$ is constructed in \cite{ARFM}. These two
    gradings should be related by Koszul duality.
  \item Part (2) of the theorem can be seen as evidence for the
    philosophy that Kazhdan-Lusztig polynomials should be replaced by
    $p$-Kazhdan-Lusztig polynomials in modular representation theory.
  \item In the analogous setting of quantum groups at a root of unity
    part (2) of the theorem (with $p$-Kazhdan-Lusztig polynomials
    replaced by ordinary Kazhdan-Lusztig polynomials) is a theorem of Soergel \cite{SoeKL,SoergelKippKM}.
  \item As we discussed in \S\ref{sec:tilting_char}, if $p \ge 2h -2$ then a small part of the knowledge of tilting
    characters can be used to obtain the simple characters. Thus the
    above theorem implies a (rather complicated) formula for the
    simple characters in terms of anti-spherical $p$-Kazhdan-Lusztig
    polynomials. It is not difficult to see that this formula implies
    Lusztig's conjecture for large $p$. However, this formula is
    \emph{not} simply the Lusztig character formula \eqref{LCF} 
    with Kazhdan-Lusztig polynomials replaced by $p$-Kazhdan-Lusztig polynomials.
  \end{enumerate}
\end{remark}


\section*{List of notation} Here is a list of frequently used
notation, in order of appearance:



\begin{longtabu}{r|l} 

$X, X^\vee$ & character lattice, cocharacter lattice, \S\ref{sec:data}
\\
$\Phi, \Phi^\vee$ & roots, coroots, \S\ref{sec:data}
\\
$G_\ZM$ & Chevalley group scheme corresponding to our root datum, \S\ref{sec:data}
\\
$\bk$, $p$ & an algebraically closed field, its characteristic, \S\ref{sec:data}
\\
$G_\bk$ & our connected, semi-simple and simply connected group over
$\bk$, \S\ref{sec:data}
\\
$T_\bk, B_\bk $ & maximal torus, Borel subgroup in $G_\bk$, \S\ref{sec:data}
\\
$\Phi_+, \Phi_+^\vee$ & positive roots, positive coroots, \S\ref{sec:data}
\\
$X_+, X_+^\vee$ & dominant weights and coweights,
\S\ref{sec:data}
\\
$\Rep H$ & abelian category of algebraic representations of $H$,
\S\ref{sec:simples}
\\
$\Irr H$ & isomorphism classes of simple $H$-modules,
\S\ref{sec:simples}
\\
$\Rep$ & algebraic representations of $G_\bk$,
\S\ref{sec:simples}
\\
$L_\l$ & simple module with highest weight $\l \in X$,
\S\ref{sec:simples}
\\
$\Delta_\l, \nabla_\l$ & Weyl and induced module with highest weight $\l \in X$,
\S\ref{sec:simples}
\\
$\DM$ & a duality on $\Rep$ fixing simples,
\S\ref{sec:simples}
\\
$W, S$ & the Weyl group and its simple reflections, \S\ref{sec:characters}
\\
$\wdot, \rho$ & the dot action, the half sum of positive roots, \S\ref{sec:characters}
\\
$\e_x$ & the sign of $x \in W$, \S\ref{sec:characters}
\\
$\ch, \chi_\l$ & the character, the Weyl character, \S\ref{sec:characters}
\\
$X_1^\ell$ & $\ell$-restricted weights, \S\ref{sec:steinberg}
\\
$(-)^\Fr$ & Frobenius twist functor, \S\ref{sec:steinberg}
\\
$T_\l$ & indecomposable tilting module,
\S\ref{sec:tilting}
\\
$\Wa, \Sa$ & affine Weyl group, its simple reflections, \S\ref{sec:Wa}
\\
$\Wae, \Omega$ & extended affine Weyl group, its length zero
elements, \S\ref{sec:Wa}
\\
$\Hee$ & extended affine Hecke algebra,
\S\ref{sec:hecke}
\\
$\He, \He_f$ & affine Hecke algebra, finite Hecke algebra,
\S\ref{sec:hecke}
\\
$h_x$ & standard basis, \S\ref{sec:hecke}
\\
$\un{h}_x, h_{y,x}$ & Kazhdan-Lusztig basis, Kazhdan-Lusztig polynomial, \S\ref{sec:hecke}
\\
$\wdot, \pdot$ & dot action, $p$-dilated dot action, \S\ref{sec:linkage}
\\
$C_-^p$ & fundamental domain for $p$-dilated dot action, \S\ref{sec:linkage}
\\
$C_+^p$ & dominant weights in the smallest $p$ alcove, \S\ref{sec:linkage}
\\
$\Rep_\mu$ & block of $\Rep$,
\S\ref{sec:linkage}
\\
$\Repp$ & principal block, \S\ref{sec:repp}
\\
$\Repe$ & extended principal block, \S\ref{sec:repp}
\\
$T_\l^\mu$ & translation functor, \S\ref{sec:translation}
\\
$h$ & Coxeter number, \S\ref{sec:translation}
\\
$D_H(X)$ & equivariant derived category, \S\ref{sec:sheafb}
\\
$D_{(H)}(X)$ & constructible derived category, \S\ref{sec:sheafb}
\\
$\Perv_{H}(X)$ & equivariant perverse sheaves,
\S\ref{sec:sheafb}\\
$\Perv_{(H)}(X)$ & perverse sheaves, \S\ref{sec:sheafb}
\\
$\Delta_\l, \Delta_\l^\bk$ & standard sheaf, \S\ref{sec:sheafb} \\
$\nabla_\l, \nabla_\l^\bk$ & costandard sheaf, \S\ref{sec:sheafb}
\\
$\ic_\l, \ic_\l^\bk$ & Intersection cohomology sheaf, \S\ref{sec:sheafb}
\\
$\EC_\l,  \EC_\l^\bk$ & parity sheaf, \S\ref{sec:sheafb}
\\
$G^\vee$ & (complex) dual group, \S\ref{sec:affineGrassFlag}
\\
$G^\vee((t))$ & loop group, \S\ref{sec:affineGrassFlag}
\\
$K,  \Iw$ & maximal compact subgroup, Iwahori subgroup, \S\ref{sec:affineGrassFlag}
\\
$\Flx,  \Grx$ & Bruhat cells, \S\ref{sec:affineGrassFlag} \\
$t^\l,  \Grl$ & special point associated to $\l \in X$, its
$K$-orbit, \S\ref{sec:affineGrassFlag} 
\\
$*$ & convolution (on affine Flag variety or affine Grassmannian),
\S\ref{sec:heckecat}, \S\ref{sec:Satake},  \S\ref{sec:FM}
\\
$\HC,\HCe$ & Hecke category, extended Hecke category, \S\ref{sec:heckecat}

\end{longtabu}

\def\cprime{$'$} \def\cprime{$'$} \def\cprime{$'$}




\end{document}